\documentclass[MSNbibl,number,citesort,seceqn,dvips]{arxbj}
\usepackage{upgreek}
\usepackage{mathbh}

\aid{0}
\volume{20}
\issue{3}
\pubyear{2014}
\firstpage{1210}
\lastpage{1233}
\doi{10.3150/13-BEJ519} 

\makeatletter
\newcommand{\RMe}{\mathrm{e}}

\newcommand{\mrmd}{\,\mathrm{d}}
\newcommand{\mrmdd}{\mathrm{d}}

\newproclaim{definition}{Definition}
\newremark{remark}[definition]{Remark}
\newremark{example}[definition]{Example}

\newtheorem{theorem}[definition]{Theorem}
\newtheorem{proposition}[definition]{Proposition}
\newtheorem{corollary}[definition]{Corollary}

\newcommand{\R}{\mathbb{R}}
\renewcommand{\SS}{\mathbb{S}}
\newcommand{\E}{{\mathbf E}}
\renewcommand{\P}{{\mathbf P}}
\newcommand{\Q}{{\mathbf Q}}
\newcommand{\one}{\mathbh1}
\newcommand{\eps}{\varepsilon}

\newcommand{\salg}{\mathfrak{K}}
\newcommand{\dsim}{\stackrel{d}{\sim}}
\newcommand{\Z}{\mathbb{Z}}
\newcommand{\kappat}{\tilde{\kappa}}
\newcommand{\imagi}{\bolds{\imath}}
\newcommand{\sH}{\mathcal{H}}
\newcommand{\sE}{\mathcal{E}}
\newcommand{\sI}{\mathcal{I}}
\newcommand{\sT}{\mathcal{T}}

\makeatother

\begin{document}
\begin{frontmatter}

\title{Invariance properties of random vectors and stochastic
processes based on the zonoid concept}
\runtitle{Invariance properties}

\begin{aug}
\author{\inits{I.}\fnms{Ilya} \snm{Molchanov}\corref{}\thanksref{e1}\ead[label=e1,mark]{ilya.molchanov@stat.unibe.ch}},
\author{\inits{M.}\fnms{Michael} \snm{Schmutz}\thanksref{e2}\ead[label=e2,mark]{michael.schmutz@stat.unibe.ch}} \and
\author{\inits{K.}\fnms{Kaspar} \snm{Stucki}\thanksref{e3}\ead[label=e3,mark]{kaspar.stucki@stat.unibe.ch}}
\runauthor{I. Molchanov, M. Schmutz and K. Stucki} 
\address{Institute of Mathematical Statistics and Actuarial Science,
University of Bern,
Sidlerstrasse 12, 3012 Bern, Switzerland. \printead{e1,e2};\\
\printead*{e3}}
\end{aug}

\received{\smonth{4} \syear{2012}}
\revised{\smonth{12} \syear{2012}}

%
\begin{abstract}
Two integrable random vectors $\xi$ and $\xi^*$ in $\R^d$ are said
to be zonoid equivalent if, for each $u\in\R^d$, the scalar products
$\langle\xi,u\rangle$ and $\langle\xi^*,u\rangle$ have the same
first absolute moments. The paper analyses stochastic processes
whose finite-dimensional distributions are zonoid equivalent with
respect to time shift (zonoid stationarity) and permutation of its
components (swap invariance). While the first concept is weaker than
the stationarity, the second one is a weakening of the
exchangeability property. It is shown that nonetheless the ergodic
theorem holds for swap-invariant sequences and the limits are
characterised.
\end{abstract}

%

\end{frontmatter}

\section{Introduction}
\label{secintroduction}

The first absolute moments $\E|\langle\xi,u\rangle|$, $u\in\R^d$, for
the scalar product of an integrable random vector $\xi$ in $\R^d$ and
$u$, admit a straightforward geometric interpretation as the support
function of a zonoid of $\xi$, see \cite{mos02}. Zonoids form an
important family of convex bodies (i.e., convex compact sets) in the
Euclidean space $\R^d$, see \cite{schn}. \emph{Zonoids} are obtained
as limits of zonotopes in the Hausdorff metric, while zonotopes are
Minkowski (elementwise) sums of a finite number of segments.

The sums of segments and the limits of sums can be interpreted as
expectations of random segments. By translation, it is possible to
assume that all segments are centred and so are of the form
$[-\xi,\xi]$ for a random vector $\xi\in\R^d$. Recall that the
\emph{support function} of a set $K$ in $\R^d$ is given by
\[
h_K(u)=\sup\bigl\{\langle u,x\rangle\dvt  x\in K\bigr\},\qquad u\in
\R^d,
\]
where $\langle u,x\rangle$ denotes the scalar product. The
expectation of $[-\xi,\xi]$ is the convex set $Z_\xi^o$ identified by
its support function, which is equal to the expected support function
of the segment (see \cite{mo1}, Section 2.1), that is,
\[
h_{Z_\xi^o}(u)=\E\bigl|\langle u,\xi\rangle\bigr|,\qquad u\in\R^d.
\]
If $\xi$ is integrable, $Z_\xi^o$ is an origin symmetric \emph{convex
body} (compact convex set). For instance, if $\xi$ is discrete in
$\R^2$ with only two possible values, then $Z_\xi^o$ is a
parallelogram; if $\xi$ is isotropic, then $Z_\xi^o$ is a ball.

A slightly different construction of zonoids associated with random
vectors was suggested by Koshevoy and Mosler, see \cite{kosmos98} and
\cite{mos02}. Namely, the zonoid $Z_\xi$ of $\xi$ is the expectation
of $[0,\xi]$ and so the support function of $Z_\xi$ is given by
\[
h_{Z_\xi}(u)=\E\langle u,\xi\rangle_+,\qquad u\in\R^d,
\]
where $x_+=\max(x,0)$. In order to stress the difference between the
two variants of zonoids, we call $Z_\xi^o$ the \emph{centred zonoid} of
$\xi$, see Section \ref{secnon-centred-zonoids} for the comparison of
the two concepts. Note that $Z_\xi$ is also well defined for some
non-integrable $\xi$. Nonetheless from now on we always assume that
all mentioned random variables and random vectors are integrable and
not identically zero.

It is well known that the zonoid of $\xi$ does not uniquely
characterise its distribution. For instance, on the line, $Z_\xi$ is
the segment with end-points determined by the expectations of the
positive and negative parts of $\xi$, while $Z_\xi^o$ is the segment
with end-points $\pm\E|\xi|$. Thus, all random variables with the same
first absolute moment are not distinguishable in terms of their
centred zonoids.

The concept of zonoid is useful in multivariate statistics to
define trimming and data depth, see \cite{cas10,mos02}.
In case of (non-centred) zonoids, the expectations
$h(k,u)=\E(k+\langle u,\xi\rangle)_+$ for $k\in\R$ and $u\in\R^d$
uniquely determine the distribution of $\xi$, and determine the
support function of a convex body in $\R^{d+1}$ called the \emph{lift
zonoid} of $\xi$, see \cite{kosmos98,mos02}.
In finance, $\E(k+\langle u,\xi\rangle)_+$ becomes the non-discounted
price of a basket call option with strike $-k$ for $k\leq0$ (if the
expectation is taken with respect to a chosen martingale measure).
The well-known result of Breeden and Litzenberger \cite{brelit78}
saying that the prices of all call options determine the distribution
of $\xi$ now becomes a corollary of a general uniqueness result for
lift zonoids, see \cite{mos02}, Theorem 2.21, and \cite{molsch10}.

\begin{definition}
\label{defze}
Two integrable random vectors $\xi$ and $\xi^*$ in $\R^d$ are called
\emph{zonoid equivalent} if their centred zonoids coincide,
that is,
\[
\E\bigl|\langle u,\xi\rangle\bigr|=\E\bigl|\bigl\langle u,\xi^*\bigr\rangle\bigr|
\]
for all $u\in\R^d$. Two families of integrable random variables
$\{\xi_t, t\in T\}$ and $\{\xi^*_t, t\in T\}$ are called
\emph{zonoid equivalent} if all their finite-dimensional
distributions are zonoid equivalent.
\end{definition}

The concept of zonoid equivalence is closely related to spectral
representations of symmetric stable ($S\alpha S$) and max-stable
processes. For instance, each $S\alpha S$ process with
$\alpha\in(0,2)$ admits the spectral representation
%
\begin{equation}
\label{eqspec-rep} X_t\dsim\int_E
f_t(z) M_\alpha(\mrmdd z),\qquad t\in T,
\end{equation}
where the equality is understood in the sense of all
finite-dimensional distributions, $\{f_t, t\in T\}$ is a family of
functions from $L^\alpha(E,\mathcal{E},\mu)$ for a measurable space
$(E,\mathcal{E},\mu)$ and $M_\alpha$ is an $S\alpha S$ random measure
with control measure $\mu$, see \cite{samtaq94}. If $X_t$ admits
another spectral representation on a measurable space
$(G,\mathcal{G},\nu)$ with functions $\{g_t\}$, then the collections
of functions $\{f_t\}$ and $\{g_t\}$ satisfy
%
\begin{equation}
\label{eqspec-chf} \int_E \Biggl|\sum
_{i=1}^n u_i f_{t_i}
\Biggr|^\alpha \mrmd \mu=\int_G \Biggl|\sum
_{i=1}^n u_i g_{t_i}
\Biggr|^\alpha \mrmd \nu
\end{equation}
for all $n\geq1$, $u_1,\ldots,u_n\in\R$ and $t_1,\ldots,t_n\in T$. This
is easily seen by computing the characteristic function of the
spectral representations, see \cite{samtaq94}, Section 3.2. If
$\alpha=1$ and both $\mu$ and $\nu$ are probability measures, then
(\ref{eqspec-chf}) can be interpreted as the zonoid equivalence of
stochastic processes $\{f_t\}$ and $\{g_t\}$.

Fairly similar facts hold for max-stable processes, see
\cite{kab09s,kabstoev12,wanstoev10,wanstoev10s}. This close
relationship between stable processes and zonoid equivalence makes it
possible to figure out a number of properties of stochastic processes
in relation to their zonoid equivalence.

The paper starts with the analysis of the main implication of the
zonoid equivalence. Namely, in Section \ref{seczono-equiv-homog} we
show that the zonoid equivalence yields the equality of the expected
values for all even one-homogeneous function of the random vectors.
Stochastic processes whose finite-dimensional distributions remain
zonoid equivalent for time shifts are discussed in
Section \ref{secisometries}. This zonoid stationarity property is
brought in relationship to the stationarity of related stable and
max-stable processes through their LePage representations.

A result of Hardin (\cite{har81}, Theorem 1.1) implies that the distribution
of an integrable random vector $\xi$ is uniquely determined by
$\E|1+\langle u,\xi\rangle|$ for all $u\in\R^d$, equivalently by the
centred zonoid of $(1,\xi)$. In Theorem \ref{corpm}, we show that, if
$\xi$ is symmetric, it is possible to replace $1$ by any random
variable taking values $\pm1$.

Section \ref{secsequences} introduces the swap-invariance property
for a random sequence that amounts to the zonoid equivalence of each
permutation of all its finite subsequences, which is weaker than the
exchangeability property. We prove the ergodic theorem for
swap-invariant sequences and characterise the limits, thereby
generalising the classical results for exchangeable sequences. Zonoid
equivalence of positive random vectors with respect to permutation of
two their components has been investigated in \cite{molsch11}
and for all possible permutations in \cite{molsch11m} in view of
financial applications.

Section \ref{secnon-centred-zonoids} discusses relationships between
centred and non-centred zonoids and also another symmetry property
being stronger than the exchangeability. In this relation, consider
\[
\E|u_0+u_1\xi_1+\cdots+u_d
\xi_d|
\]
as function $f(u_0,u_1,\ldots,u_d)$ of $(d+1)$ real arguments. The swap
invariance means exactly that $f$ is invariant for permutations of
$u_1,\ldots,u_d$ with $u_0=0$; the exchangeability corresponds to the
permutation invariance of $u_1,\ldots,u_d$ for any (and then all)
$u_0\neq0$. Assuming the full permutation invariance for all
$u_0,u_1,\ldots,u_d$ imposes a property (called \emph{lift
swap-invariance}), which is stronger than the exchangeability of
$\xi_1,\ldots,\xi_d$. A variant of this property for non-centred
zonoids has been considered in \cite{molsch10} and \cite{molsch11m}
motivated by applications in finance.

Finally, Section \ref{secequality-zonoids} collects a number of
relevant results concerning zonoids of particular distributions. It
is shown that zonoids identify uniquely distributions from
location-scale families under rather mild conditions. The special case
of random vectors with positive coordinates is also analysed, in
particular log-infinitely divisible laws being important in financial
applications.

The consideration of (non-centred) zonoids makes it possible to study
possibly non-integrable random vectors, which is left for a future
work. The same relates to $L^p$-zonoids considered in~\cite{mo09}. A
number of results of this paper can be generalised for random elements
in Banach spaces along the lines of \cite{bor09}.

\section{Expectations of homogeneous functions}
\label{seczono-equiv-homog}

Let $\sH$ (resp., $\sH_e$) denote the family of all
(resp., even) measurable homogeneous functions
\mbox{$\R^d\mapsto\R_+$}, so that $f(cx)=cf(x)$ for all $x\in\R^d$ and
$c\geq
0$.

\begin{theorem}
\label{thrhom}
Two random vectors $\xi$ and $\xi^*$ are zonoid equivalent if and
only if $\E f(\xi)=\E f(\xi^*)$ for all $f\in\sH_e$.
\end{theorem}
\begin{pf}
\textsl{Sufficiency} is immediate, since $f(x)=|\langle u,x\rangle|$
belongs to $\sH_e$.

\textsl{Necessity.} First, show that $\E\|\xi\|=\E\|\xi^*\|$. The
integral of the support function of a convex body $K$ over the unit
sphere is $\frac{1}{2} d\kappa_d b(K)$, where $b(K)$ is called the mean
width of $K$ and $\kappa_d$ is the volume of the unit ball in
$\R^d$. By changing the order of integral and expectation, it is
easy to see that the mean width of $Z^o_\xi$ equals the expected
mean width of the segment $[-\xi,\xi]$. The mean width of this
segment can be found from the Steiner formula
(\cite{schn}, Equation (4.1.1)), see also~\cite{schn}, page 210, as
$b([-\xi,\xi])=4\|\xi\|\kappa_{d-1}/(d\kappa_d)$. Thus,
$\E\|\xi\|=b(Z_\xi^o)d\kappa_d/(4\kappa_{d-1})$ is uniquely
determined by $Z_\xi^o$.

Denote the common value of $\E\|\xi\|$ and $\E\|\xi^*\|$ by $c$, and
define probability measure $\Q$ with density
\[
\frac{\mrmdd\Q}{\mrmdd\P}=\frac{\|\xi\|}{c}
\]
and another measure $\Q^*$ generated by $\xi^*$ in the same
way. Denote by $\E_\Q$ the expectation with respect to $\Q$ (and,
resp., with respect to $\Q^*$). Then for all $u\in\R^d$
\[
\frac{1}{c} \E\bigl|\langle u,\xi\rangle\bigr|=\frac{1}{c} \E\bigl|\langle u,\xi
\rangle\bigr|\one_{\{ \|\xi\|\neq0\}}= \E_{\Q} \biggl|\biggl\langle u,\frac{ \xi
}{\| \xi\|}
\biggr\rangle\biggr|\one_{\{ \|\xi\|
\neq
0\}} =\E_{\Q} \biggl|\biggl\langle u,
\frac{ \xi}{\| \xi\|}\biggr\rangle\biggr|
\]
and similarly $c^{-1}\E|\langle u,\xi^*\rangle|=\E_{\Q^*}|\langle u,
\xi^*/\|\xi^*\|\rangle|$. Therefore, $\xi/\|\xi\|$ under $\Q$ and
$\xi^*/\|\xi^*\|$ under $\Q^*$ share the same zonoid. Define
measure $\mu$ on the unit Euclidean sphere by setting
$\mu(A)=\Q(\xi/\| \xi\|\in A)$ and correspondingly $\mu^*$. The
convex body $Z_\mu^o$ with the support function
\[
h_{Z_\mu^o}(u)=\int_{\SS^{d-1}} \bigl|\langle u,x\rangle\bigr| \mu(\mrmdd x) =
\E_{\Q} \biggl|\biggl\langle u,\frac{ \xi}{\| \xi\|}\biggr\rangle\biggr|
\]
is termed the zonoid of $\mu$, see \cite{schn}, Section 3.5. It is well
known that an even finite measure on the unit sphere is uniquely
determined by its zonoid, see \cite{schn}, Theorem 3.5.3. Since $\mu$
and $\mu^*$ share the same zonoid, the integrals of any even and
integrable function with respect to them coincide.

For $f\in\sH_e$, we have $f(0)=0$ and therefore
\[
\E f(\xi)=\E\bigl[f(\xi)\one_{\|\xi\|\neq0}\bigr] =\E_\Q f\bigl(\xi/\|
\xi\|\bigr)=\int_{\SS^{d-1}} f(u)\mu(\mrmdd u).
\]
Hence, $\E f(\xi)=\E f(\xi^*)$ for each $f\in\sH_e$. A short
calculation shows that integrability of $f(\xi/\| \xi\|)$ under
$\Q$ implies integrability of $f(\xi^*/\| \xi^*\|)$ under $\Q^*$ and
vice versa.
\end{pf}

If $\xi$ and $\xi^*$ are zonoid equivalent, then $f(\xi)$ and
$f(\xi^*)$ are two zonoid equivalent random variables for all
$f\in\sH_e$. The following result is easily derived by observing that
$\E f(\xi)=\E\frac{1}{2}(f(\xi)+f(-\xi))$ for symmetric $\xi$.

\begin{corollary}
\label{thrhom-sym}
Two symmetric random vectors $\xi$ and $\xi^*$ are zonoid equivalent
if and only if $\E f(\xi)=\E f(\xi^*)$ for all $f\in\sH$. In
particular, $\E h_K(\xi)=\E h_K(\xi^*)$ for each convex body~$K$.
\end{corollary}

\begin{corollary}
\label{cormany-hom}
Let $f_1,\ldots,f_k\in\sH_e$. If $\xi$ and $\xi^*$ are zonoid
equivalent, then the vectors $(f_1(\xi),\ldots,f_k(\xi))$ and
$(f_1(\xi^*),\ldots,f_k(\xi^*))$ are zonoid equivalent as long as one
of these vectors is integrable.
\end{corollary}
\begin{pf}
It suffices to use the fact that $f(x)=|u_1f_1(x)+\cdots+u_kf_k(x)|$
belongs to $\sH_e$ and $f(\xi)$ is integrable.
\end{pf}

The following easy fact is also worth noticing.

\begin{proposition}
\label{proplinear}
Two random vectors are zonoid equivalent if and only if all their
linear transformations are zonoid equivalent.
\end{proposition}
\begin{pf}
For each matrix $A$, we have $\E|\langle A\xi,u\rangle|=\E|\langle
\xi,A^\top u\rangle|$ and so $A\xi$ and $A\xi^*$ are zonoid
equivalent if $\xi$ and $\xi^*$ are.
\end{pf}

In the following, we often consider random vectors with positive
coordinates (shortly called positive vectors), which are usually
denoted by the letter $\eta$.

\begin{proposition}
\label{prophf-max}
Two positive integrable random vectors $\eta$ and $\eta^*$ are
zonoid equivalent if and only if $\E f(\eta)=\E f(\eta^*)$ for each
$f\in\sH$. In particular, the zonoid equivalence implies
$\E\eta=\E\eta^*$.
\end{proposition}
\begin{pf}
While the sufficiency is evident, the necessity can be proved
similarly to Theorem \ref{thrhom} with $\Q$ having density
$\eta_1/\E\eta_1$. The equality of expectations is obtained by
setting $f(x)=(x_i)_+$ for any $i=1,\ldots,d$.
\end{pf}

For positive random vectors, the concept of a \emph{max-zonoid} is
also useful. The max-zonoid $M_\eta$ of a positive random vector
$\eta=(\eta_1,\ldots,\eta_d)$ is defined as the expectation of the
crosspolytope in $\R^d$, which is the convex hull of the origin and
the standard basis vectors scaled by $\eta_1,\ldots,\eta_d$, see
\cite{mo08e}. The support function of $M_\eta$ is given by
%
\begin{equation}
\label{eqsup-max-z} h_{M_\eta}(u)=\E\max(0,u_1
\eta_1,\ldots,u_d\eta_d),\qquad u=(u_1,\ldots,u_d)\in\R^d.
\end{equation}
This support function is most interesting for positive
$u_1,\ldots,u_d$, where it is possible to omit zero in the right-hand
side of (\ref{eqsup-max-z}). The following result has been proved
analytically in \cite{wanstoev10}, Theorem 1.1. An alternative proof
(using a geometric argument combined with the change of measure
technique) has recently been given in \cite{molsch11m}, Proposition 1.

\begin{proposition}
\label{propmz-eq}
Two positive integrable random vectors $\eta$ and $\eta^*$ have
identical max-zonoids if and only if $\eta$ and $\eta^*$ are zonoid
equivalent.
\end{proposition}

\section{Isometries, representations of stable processes, and zonoid
stationarity}
\label{secisometries}

A result of Hardin (\cite{har81}, Theorem 1.1) reformulated for random
vectors implies that, for any given positive $p\notin2\Z$, the values
$\E|1+\langle u,\xi\rangle|^p$ for all $u\in\R^d$ determine uniquely
the distribution of random vector $\xi\in\R^d$. If $p=1$, this result
means that the centred zonoid of $(1,\xi)$ uniquely identifies the
distribution of $\xi$, cf. \cite{kosmos98,mos02}. This also means
that if two zonoid equivalent random vectors contain the same
coordinate being exactly one, then these random vectors are
identically distributed. Below we provide a generalisation of this
result for $p=1$ and symmetric random vectors showing that it is
possible to replace the constant $1$ with any random variable taking
values $\pm1$.

\begin{theorem}
\label{corpm}
Let $\xi$ be a symmetric random vector in $\R^d$. If $\eps$ is any
random variable with values~$\pm1$, then the centred zonoid of
$(\eps,\xi)$, that is, the values of
\[
\E\bigl|u_0\eps+\langle u,\xi\rangle\bigr|,\qquad u_0\in\R, u\in
\R^d,
\]
determines uniquely the distribution of $\xi$.
\end{theorem}
\begin{pf}
For each function $f(\eps,\xi)$ we have
$f(\eps,\xi)+f(-\eps,\xi)=f(1,\xi)+f(-1,\xi)$, so that
\[
\E\bigl|u_0\eps+\langle u,\xi\rangle\bigr|+\E\bigl|-u_0\eps+\langle
u,\xi\rangle\bigr| =\E\bigl|u_0+\langle u,\xi\rangle\bigr|+\E\bigl|-u_0+
\langle u,\xi\rangle\bigr|.
\]
Since $\xi$ is symmetric,
\[
\E\bigl|-u_0+\langle u,\xi\rangle\bigr|=\E\bigl|u_0+\langle u,-\xi
\rangle\bigr| =\E\bigl|u_0+\langle u,\xi\rangle\bigr|.
\]
Thus,
\[
\E\bigl|u_0+\langle u,\xi\rangle\bigr| =\tfrac{1}{2}\bigl(
\E\bigl|u_0\eps+\langle u,\xi\rangle\bigr| +\E\bigl|-u_0\eps+\langle u,
\xi\rangle\bigr|\bigr)
\]
for all $u_0\neq0$ and $u\in\R^d$. Therefore, the right-hand side
is determined by the zonoid of $(\eps,\xi)$, and it remains to note
that the left-hand side uniquely identifies the distribution of
$\xi$ by \cite{har81}, Theorem 1.1.
\end{pf}

An integrable random vector $\xi$ in $\R^d$, which is not a.s.
zero, generates a norm on $\R^d$ by
\[
\|u\|_\xi=\E\bigl|\langle u,\xi\rangle\bigr|.
\]
With this definition, zonoid equivalence of $\xi$ and $\xi^*$ means
that \mbox{$\|\cdot\|_\xi$} and \mbox{$\|\cdot\|_{\xi^*}$} are two identical norms
on $\R^d$. The uniqueness result in \cite{har81} is used to
characterise isometries of subspaces of $L^1$ that contain the
function identically equal one. Theorem \ref{corpm} makes it possible
to obtain similar results for subspaces of $L^1$ that consist of
symmetric random variables and contain random variables taking values
$\pm c$ for any fixed $c>0$. The characterisation of linear
isometries defined on families of random variables are important for
the studies of symmetric stable laws, see \mbox{\cite{har81,har82,ros06}}.

A collection of integrable random elements $\{\xi_t, t\in T\}$ is a
subset of the space $L^1=L^1(\Omega,\salg,\P)$. Denote by $F_\xi$ the
$L^1$-closure of the linear space generated by this collection. Assume
that $\Omega$ is a Borel space with $\salg$ being the Borel
$\sigma$-algebra.

Assume that $\{\xi_t\}$ is rigid, that is, any linear isometry
$U_0\dvtx
F_\xi\mapsto L^1$ is uniquely extendable to the isometry $U\dvtx
L^1\mapsto L^1$. It is well known \cite{har81,ros06} that the rigidity
is guaranteed by imposing that the random elements $\{\xi_t\}$ have
full support, the union of its supports is $\Omega$ up to a null set
(see \cite{har82} for details), and that $\xi_t/\bar{\xi}$, $t\in T$,
generate the $\sigma$-algebra $\salg$, where $\bar{\xi}\in F_\xi
$ is
a random variable with full support (its existence is guaranteed by
\cite{har81}, Lemma 3.2). Note that the family $\{\xi_t\}$ is often
called \emph{minimal} instead of \emph{rigid}, as it gives rise to a
\emph{minimal spectral representation} of a $S\alpha S$ process via
(\ref{eqspec-rep}), see also \cite{har82}.

Consider another rigid collection $\{\xi^*_t, t\in T\}$, which is
zonoid equivalent to $\{\xi_t, t\in T\}$. Then the
isometry between $F_\xi$ and $F_{\xi^*}$ can be characterised as
follows, see Theorem 3.2 in \cite{ros06}. For every $t\in T$,
%
\begin{equation}
\label{eqiso} \xi^*_t(\omega)=h(\omega)\xi_t\bigl(
\phi(\omega)\bigr) \qquad\P\mbox{-a.s.,}
\end{equation}
where $\phi\dvtx \Omega\to\Omega$, $h\dvtx \Omega\to\R\setminus\{0\}$ are
measurable and $|\bar{\xi}|\mrmdd\P=|\bar{\xi}|(|h|\mrmdd\P) \circ\phi^{-1}$,
for a random variable  $\bar{\xi}\in F_\xi$ with full support.

A similar construction of isometries can be carried over for
max-zonoids and non-negative integrable functions, see
\cite{haapic86}, where such isometries are called pistons. Since for
positive random vectors the zonoid equivalence and the max-zonoid
equivalence are identical (see Proposition \ref{propmz-eq}), the
isometries corresponding to max-zonoids are also characterised by
(\ref{eqiso}).

Recall that each \emph{symmetric $1$-stable} (i.e., $S\alpha S$ with
$\alpha=1$) random vector $X$ in $\R^d$ can be represented as the
LePage series
%
\begin{equation}
\label{eq3} X=\sum_{k=1}^\infty
\Gamma_k^{-1} \xi^{(k)},
\end{equation}
where $\Gamma_k=\zeta_1+\cdots+\zeta_k$ are successive sums of
i.i.d.
standard exponential random variables and
$\xi,\xi^{(1)},\xi^{(2)},\ldots$ are i.i.d. integrable symmetric
random vectors, see \cite{lepwoodzin81}. Note that the $\xi$'s are
often assumed to be distributed on the unit sphere with an extra
normalisation constant in front of the sum, see
\cite{samtaq94}, Corollary 1.4.3. A similar series representation with the
sum replaced by the maximum, and positive $\xi$ yields simple
(i.e., having unit Fr\'echet marginals) max-stable random vectors, see
\cite{haa84} and \cite{falhusrei04}. If $\xi$ is a stochastic
process, similar series representations yield symmetric $1$-stable
processes and simple max-stable processes. For instance, a
result of \cite{haa84} says that each stochastically continuous simple
max-stable process $Y$ can be represented as
%
\begin{equation}
\label{eqBR} Y_t=\max_{k\geq1} \Gamma_k^{-1}
\xi^{(k)}_t,\qquad t \in\R^d,
\end{equation}
where $\{\xi^{(k)}_t, t\in\R^d\}$ are i.i.d. copies of an
integrable positive process $\{\xi_t, t\in\R^d\}$. In the following,
we refer to (\ref{eq3}), its variant for stochastic processes or
their max-analogues as the LePage series.

\begin{theorem}
\label{thrlepage}
Two LePage series $X$ and $X^*$ given by (\ref{eq3}) (resp.,
their max-analogues) with integrable symmetric (resp.,
positive) summands distributed as $\xi$ and $\xi^*$ coincide in
distribution if and only if $\xi$ and $\xi^*$ are zonoid equivalent.
\end{theorem}
\begin{pf}
It suffices to consider the case of $\xi$ being a random vector in
$\R^d$. The points $\{(\Gamma_k^{-1},\xi^{(k)}), k\geq1\}$ build
the Poisson point process on $(0,\infty)$ with intensity $t^{-2}$,
$t>0$, and independent marks $\xi^{(k)}$, $k\geq1$. The formula for
the probability generating functional of the marked Poisson process
(see \cite{davj}) yields the characteristic function of $X$
\begin{eqnarray*}
\E \RMe^{\imagi\langle u,X\rangle} &=&
\exp\biggl\{-\int_0^\infty
\E\bigl(1-\RMe^{\imagi t\langle u,\xi\rangle}\bigr) t^{-2}\mrmd t\biggr\}
\\
&=&
\exp\biggl\{-\int_0^\infty\E\bigl(1-\cos\bigl(t
\langle u,\xi\rangle\bigr)\bigr) t^{-2}\mrmd t\biggr\}=\exp\biggl\{-
\frac{\uppi}{2}\E\bigl|\langle u,\xi\rangle\bigr|\biggr\},
\end{eqnarray*}
since $\int_0^\infty(1-\cos(s))s^{-2}\mrmd s=\uppi/2$, where $\imagi$
denotes the imaginary unit. Thus, the distribution of $X$ is
determined by $\E|\langle u,\xi\rangle|$, $u\in\R^d$.

The result for max-stable random vectors follows from the
association argument from \cite{kab09s} or \cite{wanstoev10} or a
direct calculation of the cumulative distribution functions combined
with Proposition~\ref{propmz-eq}.
\end{pf}

Let $\{\xi_t, t\in T\}$ be a stochastic process such that $\xi_t$ is
integrable for all $t \in T$, where $T$ is either integer grid $\Z^d$
or $\R^d$.

\begin{definition}
\label{defzstat}
The process $\{\xi_t, t\in T\}$ is called \emph{zonoid stationary}
if $\{\xi_t, t\in T\}$ and $\{\xi_{t+s}, t\in T\}$ are zonoid
equivalent for all $s \in T$.
\end{definition}

Obviously all integrable stationary processes are zonoid stationary.
If both $\xi$ and $\xi^*$ are centred Gaussian processes, then by
Corollary \ref{cormultivariate-normal} their zonoid equivalence
implies the equality of all finite-dimensional distributions, so their
zonoid stationarity is equivalent to the conventional
stationarity. The same holds for symmetric $\alpha$-stable processes
with given $\alpha>1$. The fact that zonoid does not uniquely
determine the general distribution suggests that there
exist non-stationary but zonoid stationary processes. The next
result follows from Theorem \ref{thrlepage}.

\begin{corollary}
\label{corz-st}
A symmetric 1-stable process (resp., max-stable process with
unit Fr\'echet marginals) obtained as the LePage series (\ref{eq3})
(resp., (\ref{eqBR})) is stationary if and only if $\xi$ is
zonoid stationary.
\end{corollary}

If the max-stable process $Y$ given by (\ref{eqBR}) is stationary,
the process $\log\xi$ is called \emph{Brown--Resnick stationary},
see \cite{kabschhaan09}. Corollary \ref{corz-st} shows that the
Brown--Resnick stationarity of $\log\xi$ is equivalent to the zonoid
stationarity of a positive stochastic process $\xi$.

\begin{example}
\label{exgbm}
The geometric Brownian motion $\RMe^{W_t-|t|/2}$, where $W_t$,
$t\in\R$, is a double-sided Brownian motion, is zonoid
stationary. The corresponding stationary process given by
(\ref{eqBR}) was introduced by Brown and Resnick
\cite{brores77}. Kabluchko \textit{et al.}  \cite{kabschhaan09} replaced
$W_t$ by a Gaussian process $\xi_t$ with mean $\mu_t$ and variance
$\sigma_t^2$. Their result implies that $\RMe^{\xi_t}$ is zonoid
stationary if and only if $\xi_t-\mu_t$ has stationary increments
and $\mu_t+\frac{1}{2}\sigma_t^2$ is constant for all $t$.\looseness=-1
\end{example}

For a zonoid stationary process $\xi$ the spaces generated by
$\{\xi_t,t\in T\}$ and $\{\xi_{t+h},t\in T\}$ are isometric for all $h
\in T$. This gives rise to a representation of $\xi$ in term of
isometries. Following \cite{ros95}, a measurable function
$\phi\dvtx \Omega\times T \to\Omega$ is said to be a measurable flow if
$\phi_{t_1+t_2}(\omega)=\phi_{t_1}(\phi_{t_2}(\omega))$ and
$\phi_0(\omega)= \omega$ for all $t_1,t_2\in T$ and $\omega\in
\Omega$. The flow $\phi$ is said to be non-singular if $\P\circ
\phi_t^{-1}$ is equivalent to $\P$ for all $t\in T$. A measurable
function $r\dvtx
\Omega\times T \to\R$ is said to be a cocycle for a measurable flow
$\phi$ if
$r_{t_1+t_2}(\omega)=r_{t_1}(\omega)r_{t_2}(\phi_{t_1}(\omega))$ for
all $t_1,t_2 \in T$ and for $\P$-almost all $\omega\in\Omega$. By
replicating the proofs of \cite{ros95}, Theorem 3.1, and
\cite{ros00}, Theorem 2.2, it is easy to show that a zonoid stationary
process $\xi$ with rigid (minimal) family $F_\xi$ satisfies
\[
\xi_t(\omega)=r_t(\omega) \biggl(\frac{\mrmdd\P\circ\phi_t}{\mrmdd\P
}
\biggr) (\omega) (\xi_0\circ\phi_t) (\omega) \qquad\P
\mbox{-a.s.},
\]
where $\{\phi_t, t\in T\}$ is a measurable non-singular flow and
$\{r_t, t\in T\}$ is a cocycle for $\phi$ taking values in
$\{-1,1\}$.

\section{Swap invariant sequences}
\label{secsequences}

A finite or infinite random sequence $\xi=(\xi_1,\xi_2,\ldots)$ of
random elements is said to be \emph{exchangeable} if its distribution
is invariant under finite permutations, that is, the distribution of
any finite subsequence is invariant under any permutation of its
elements, see, for example, \cite{kal05s}, Section~1.1.

\begin{definition}
An integrable random vector is called \emph{swap-invariant} if all
random vectors obtained by permutations of its coordinates are
zonoid equivalent. A sequence of integrable random variables is
called swap-invariant if all its finite subsequences are
swap-invariant.
\end{definition}

An integrable random vector $\xi$ with positive components exhibiting
the swap-invariance property restricted to permutation of its two
components $\xi_i$ and $\xi_j$ is called $ij$-swap-invariant. This
weaker variant of the swap-invariance property has been already
introduced and applied in a financial context in \cite{molsch11} and
\cite{schm10qf}. The swap-invariance property of the vector of asset
prices ensures that different financial derivatives share the same
price and can be freely exchanged, which is an essential tool for
semi-static hedging of barrier options, see \cite{carellgup98}.

The swap-invariance property of $\xi$ immediately implies that
$\E|\xi_1|=\cdots=\E|\xi_d|$. It is obvious that the exchangeable
sequence is swap-invariant. The following examples show that the
swap-invariance is weaker than the exchangeability property.

\begin{example}[(See \cite{dac82})]
\label{exdac}
On the probability space $\Omega=[0,1]$ with the Lebesgue measure define
%
\begin{equation}
\label{eqdce} \xi_n=n(n+1)\one_{\omega\in((n+1)^{-1},n^{-1}]},\qquad n\geq1.
\end{equation}
By a direct computation it is easy to see that
\[
\E|u_1\xi_1+\cdots+u_n\xi_n|=
\sum_{i=1}^n |u_i|,
\]
so that the sequence is indeed swap-invariant, but not
exchangeable. Further examples of this type can be obtained for
general sequences of non-negative random variables with equal
expectations and disjoint supports.
\end{example}

\begin{example}
\label{exlognsi}
Let $Z_1,Z_2,\ldots$ be a sequence of i.i.d. standard normal random
variables and let $\{b_k, k\geq1\}$ be a sequence of
real numbers such that $\sum b_k^2<\infty$. Define $\eta_i=\RMe^{\xi_i}$,
$i\geq1$, where
\[
\xi_i=Z_i+\sum_{k=1}^\infty
b_kZ_k+\mu_i
\]
and
\[
\mu_i=-\frac{1}{2}\operatorname{Var}(\xi_i)=-
\frac{1}{2}\Biggl(1+\sum_{k=1}^\infty
b_k^2+2b_i\Biggr).
\]

By Corollary \ref{corlognorm}, $\eta$ is swap-invariant. Note that
no two components $\eta_i$ and $\eta_j$ are identically distributed
unless $b_i=b_j$.
\end{example}

If the extended sequence $(1,\xi)$ (or $(\eps,\xi)$ with
$\eps\in\{-1,1\}$ and symmetric $\xi$) is swap-invariant, then $\xi$
is exchangeable. Actually, the swap invariance of such extended
sequence is stronger than the exchangeability of $\xi$, see
Section \ref{secnon-centred-zonoids}.

It is well known that each exchangeable sequence of integrable random
variables satisfies several ergodic theorems. Given an infinite random
sequence $\{\xi_n, n\geq1\}$, denote the corresponding tail
$\sigma$-algebra by $\sT_\xi$, the shift-invariant $\sigma$-algebra by
$\sI_\xi$, and the permutation-invariant $\sigma$-field by
$\sE_\xi$. These $\sigma$-algebras are identical modulo null sets
for exchangeable sequences, see \cite{kal05s}, Corollary 1.6. Since an
infinite exchangeable sequence is stationary, the following result is
a direct consequence of \cite{kalle}, Theorem 10.6,
and \cite{kal05s}, Corollary 1.6.

\begin{theorem}
\label{thmerg-kalle}
Let $\xi_1,\xi_2,\ldots$ be an exchangeable sequence of integrable
random variables. Then
\[
n^{-1}\sum_{i=1}^n
\xi_i \to\E(\xi_1\mid\sE_\xi) \qquad\mbox{a.s. and
in } L^1 \mbox{ as } n \to\infty.
\]
\end{theorem}

In the following, we extend this fact to swap-invariant
sequences. Recall that these sequences by definition consist of
integrable random variables.

\begin{theorem}
\label{thrsi-slln}
Let $\xi_1,\xi_2,\ldots$ be a swap-invariant sequence of random
variables. Then $n^{-1}(\xi_1+\cdots+\xi_n)$ converges almost surely
to an integrable random variable $X$ as $n\to\infty$.
\end{theorem}
\begin{pf}
Assume first that all random variables $\xi_1,\xi_2,\ldots$ are
symmetric and that at least one random variable (say $\xi_1$) is
non-zero with probability one. Recall that $\E|\xi_i|$ is the same
for all $i$. Define an equivalent to $\P$ probability measure $\P^1$
by
%
\begin{equation}
\label{eqq-def} \frac{\mrmdd\P^1}{\mrmdd\P}=\frac{|\xi_1|}{\E|\xi_1|}.
\end{equation}
For any finite subsequence $\xi=(\xi_1,\xi_{k_1},\ldots,\xi_{k_d})$,
%
\begin{equation}
\label{eq5} \frac{\E|\langle u,\xi\rangle|}{\E|\xi_1|} =\E_{\P^1}\biggl|u_1
\eps+u_2\frac{\xi_{k_1}}{|\xi_1|} +\cdots+u_d\frac{\xi_{k_d}}{|\xi_1|}\biggr|,
\end{equation}
where $\eps=\xi_1/|\xi_1|$ is the sign of $\xi_1$ and $\E_{\P^1}$
denotes the expectation with respect to $\P^1$. By
Theorem \ref{corpm}, the right-hand side of (\ref{eq5})
determines the distribution of $(\xi_{k_1},\ldots,\xi_{k_d})/|\xi_1|$
under~$\P^1$. By writing (\ref{eq5}) for a permutation
$\xi_{k_{i_1}},\ldots,\xi_{k_{i_d}}$ we arrive at the conclusion that the
sequence $\frac{\xi_2}{|\xi_1|},\frac{\xi_3}{|\xi_1|},\ldots$ is
exchangeable under $\P^1$. Theorem \ref{thmerg-kalle} yields that
\[
\frac{1}{n} \biggl(\frac{\xi_2}{|\xi_1|}+\cdots+\frac{\xi_n}{|\xi
_1|} \biggr)
\to Z \qquad\mbox{$\P^1$-a.s. as } n\to\infty
\]
for some random variable $Z$. Since $\P^1$ and $\P$ are equivalent,
the same holds $\P$-a.s. Thus,
\[
\frac{\xi_2+\cdots+\xi_n}{n}\to X=|\xi_1|Z \qquad\mbox{a.s. as } n\to\infty.
\]
It is obviously possible to add $\xi_1$ in the numerator without
altering the limit.

If the sequence $\{\xi_n\}$ is no longer symmetric, consider an
independent symmetric random variable $\eps$ with values $\pm
1$. Then the sequence $\{\eps\xi_n, n\geq1\}$ is symmetric and
swap-invariant, which is seen by the total probability formula. As
shown above, $\{\eps\xi_n\}$ satisfies the ergodic theorem with
limit $X_\eps$. Then the original sequence $\{\xi_n\}$ satisfies the
ergodic theorem with the limit $\eps X_\eps$ (note that $\eps$ and
$X_\eps$ may be dependent).

It remains to consider the case when all $\xi_i$ have an atom at
zero. Fix any $k \ge1$ and define a new measure $\P^k$ by
\[
\frac{\mrmdd\P^k}{\mrmdd\P}=\frac{|\xi_k|}{\E|\xi_k|}.
\]
The function $(x_1,\ldots,x_d)\mapsto
|u_1x_1+\cdots+u_dx_d|\one_{x_k\neq0}$ is in $\sH_e$, hence
\begin{eqnarray*}
&&\E|u_1\xi_1+u_2\xi_2+\cdots+u_k\xi_k+\cdots+ u_d\xi_d|
\one_{\xi_k\neq0}
\\
&&\quad=\E|u_1\xi_1+u_2\xi_{i_2}+\cdots+u_k \xi_k +\cdots+u_n\xi
_{i_d}|\one_{\xi_k\neq0}
\end{eqnarray*}
for all $u_1,\ldots,u_d \in\R$ and all permutations $i_1,\ldots,i_d$
with $i_k=k$. Thus, the sequence
$(\xi_1,\ldots,\xi_{k-1},\xi_{k+1},\ldots)/|\xi_k|$ is exchangeable
under $\P^k$. Since $\P^k$ is equivalent to $\P$ restricted on
${\{\xi_k\neq0\}}$, $n^{-1}(\xi_1+\cdots+\xi_n)$ converges to some
random variable $X$ for almost all $\omega\in\{\xi_k\neq0\}$. Note
that the same limit appears under $\P^m$ for $m\neq k$ for almost
all $\omega$ such that $\xi_k(\omega)\neq0$ and
$\xi_m(\omega)\neq0$. Finally, set $X(\omega)=0$ for all $\omega
\in
\Omega$ such that $\xi_n(\omega)=0$ for all $n\geq1$.

Since $\xi_1,\xi_2,\ldots$ have the same first absolute moment, the
integrability of $X$ follows trivially by Fatou's lemma and the
triangle inequality.
\end{pf}

\begin{remark}
A proof of Theorem \ref{thrsi-slln} for almost surely positive
swap-invariant sequences can be alternatively carried over by
using $\xi_1$ to change the measure and then referring to
\cite{har81}, Theorem~1.1.
\end{remark}

\begin{theorem}
\label{thrL1}
Assume that a swap-invariant sequence $\xi_1,\xi_2,\ldots$ satisfies
one of the following conditions:
\begin{itemize}[(b)]
\item[(a)] $\xi_k\neq0$ a.s. for some $k\ge1$,
\item[(b)] $\xi_1,\xi_2,\ldots$ is uniformly integrable.
\end{itemize}
Then the convergence of $n^{-1}(\xi_1+\cdots+\xi_n) \to X$ also holds
in $L^1$.
\end{theorem}
\begin{pf}
(a) The proofs of Theorems \ref{thrsi-slln} and
\ref{thmerg-kalle} yield that
\[
\E\bigl|n^{-1}(\xi_1+\cdots+\xi_n) - X\bigr|
\one_{\xi_k\neq0} \to0 \qquad\mbox{as } n \to\infty,
\]
while $\P(\xi_k\neq0)=1$.

(b) It is well known that the uniform integrability of
$\{\xi_n, n\ge1\}$ implies the uniform integrability of
$\{(\xi_1+\cdots+\xi_n)/n, n\ge1\}$. The a.s. convergence implies
 the $L^1$-convergence in view of
the uniform integrability property, see \cite{kalle}, Proposition 4.12.
\end{pf}

\begin{example}[(Example \ref{exdac} continuation)]
For the sequence (\ref{eqdce}), $n^{-1}(\xi_1+\cdots+\xi_n) \to0$
a.s., but $\E n^{-1}(\xi_1+\cdots+\xi_n)=1$, so the ergodic theorem
holds almost surely but not in $L^1$.
\end{example}

The following theorem characterises the limits in
Theorem \ref{thrsi-slln} for the case when at least one random
variable in the sequence does not have an atom at zero.

\begin{theorem}
\label{thmlimit}
Let $\xi=(\xi_1,\xi_2,\ldots)$ be a symmetric swap-invariant sequence
such that $\xi_1\neq0$ a.s. Then
%
\begin{equation}
\label{eqlimit} \frac{1}{n}\sum_{i=1}^n
\xi_i \to\frac{|\xi_1|}{\E(|\xi_1|\mid
\sE_{\tilde{\xi}})}\E(\xi_2\mid
\sE_{\tilde{\xi}}) \qquad\mbox{a.s. and in } L^1 \mbox{ as } n \to\infty,
\end{equation}
where $\tilde{\xi}=(\xi_2/|\xi_1|,\xi_3/|\xi_1|,\ldots)$.
\end{theorem}
\begin{pf}
The sequence $\tilde{\xi}$ is exchangeable under $\P^1$ defined by
(\ref{eqq-def}) and Theorem \ref{thmerg-kalle} implies
%
\begin{equation}
\label{eqQlimit} \frac{1}{n}\sum_{i=1}^n
\frac{\xi_i}{|\xi_1|} \to\E_{\P^1} \biggl[\frac{\xi_2}{|\xi_1|}\Bigm|
\sE_{\tilde{\xi}} \biggr] \qquad\mbox{a.s. and in } L^1 \mbox{ as } n \to
\infty.
\end{equation}
Let $Z$ be a $\sE_{\tilde{\xi}}$ measurable and
$\P^1$-integrable random variable. Then
%
\begin{equation}
\label{eqcondExp} \E_{\P^1} Z=\E\frac{|\xi_1| Z}{\E|\xi_1|}= \E\biggl
[\E\biggl(
\frac{|\xi_1|Z}{\E|\xi_1|} \Bigm| \sE_{\tilde
{\xi}} \biggr) \biggr] = \E\biggl[Z
\frac{\E(|\xi_1|\mid\sE_{\tilde{\xi}})}{\E|\xi
_1|} \biggr].
\end{equation}
Let $A \in\sE_{\tilde{\xi}}$. By the definition of the conditional
expectation
\begin{eqnarray*}
\E_{\P^1} \biggl(\one_A \E_{\P^1} \biggl(
\frac{\xi_2}{|\xi_1|}\Bigm| \sE_{\tilde{\xi}} \biggr) \biggr) &=& \E_{\P^1}\bigl(
\one_A \xi_2/\E|\xi_1|\bigr) = \E_{\P^1}
\bigl(\one_A\E\bigl(\xi_2/\E|\xi_1|\mid
\sE_{\tilde{\xi}}\bigr)\bigr)
\\
&=&\E_{\P^1} \biggl[\one_A\frac{\E(\xi_2\mid\sE_{\tilde{\xi
}})}{\E(|\xi
_1| \mid
\sE_{\tilde{\xi}})}
\frac{\E(|\xi_1|\mid\sE_{\tilde{\xi}})}{\E|\xi_1|} \biggr] = \E_{\P^1}
\biggl[ \one_A
\frac{\E(
\xi_2\mid\sE_{\tilde{\xi}})}{\E(|\xi_1|\mid\sE_{\tilde{\xi}})} \biggr],
\end{eqnarray*}
where the last equality follows from (\ref{eqcondExp}). The
uniqueness of the conditional expectation yields
\[
\E_{\P^1} \biggl[\frac{\xi_2}{|\xi_1|}\Bigm| \sE_{\tilde{\xi}} \biggr] =
\frac{\E(\xi_2\mid\sE_{\tilde{\xi}})}{\E(|\xi_1|\mid\sE
_{\tilde{\xi}})} \qquad\mbox{a.s.}
\]
This equation together with (\ref{eqQlimit})
yield the claim.
\end{pf}

With a similar proof, we arrive at the following result for positive
sequences.

\begin{proposition}
\label{proplimit-pos}
Let $\xi=(\xi_1,\xi_2,\ldots)$ be a positive swap-invariant sequence.
Then
%
\begin{equation}
\label{eqlimit-pos} \frac{1}{n}\sum_{i=1}^n
\xi_i \to\frac{\xi_1}{\E(\xi_1\mid
\sE_{\tilde{\xi}})}\E(\xi_2\mid
\sE_{\tilde{\xi}}) \qquad\mbox{a.s. and in } L^1 \mbox{ as } n \to\infty,
\end{equation}
where $\tilde{\xi}=(\xi_2/\xi_1,\xi_3/\xi_1,\ldots)$.
\end{proposition}

For non-symmetric swap-invariant sequences, we obtain the following
result by applying the total probability formula and
Theorem \ref{thmlimit}.

\begin{corollary}
\label{corlimit-general}
Let $\xi=(\xi_1,\xi_2,\ldots)$ be a swap-invariant sequence
such that $\xi_1\neq0$ a.s. Then
%
\begin{equation}
\label{eqlimit-gen} \frac{1}{n}\sum_{i=1}^n
\xi_i \to\frac{|\xi_1|}{\E(|\xi_1|\mid
\sE_{{\eps\tilde\xi}})\eps}\E(\eps\xi_2\mid
\sE_{\eps\tilde{\xi}}) \qquad\mbox{a.s. and in } L^1 \mbox{ as } n \to\infty,
\end{equation}
where $\eps$ is the Rademacher random variable independent of $\xi$
under $\P$.
\end{corollary}

\begin{corollary}
Let $\xi=(\xi_1,\xi_2,\ldots)$ be a swap-invariant sequence. If
$n^{-1}(\xi_1+\cdots+\xi_n)$ converges in $L^1$ to a deterministic
non-zero limit $c$, then $(c,\xi)$ is swap-invariant and so $\xi$ is
exchangeable.
\end{corollary}
\begin{pf}
For $m,n\ge1$, the swap-invariance property implies
\[
\E\Biggl|u_1\xi_1+\cdots+u_n\xi_n+u_0
\frac{1}{m}\sum_{k=1}^m
\xi_{n+k}\Biggr| = \E\Biggl|u_{i_1}\xi_1+\cdots+u_{i_n}
\xi_n+u_{i_0}\frac{1}{m}\sum
_{k=1}^m\xi_{n+k}\Biggr|
\]
for all permutations $(i_0,i_1,\ldots,i_n)$ of $(0,1,\ldots,n)$. The
$L^1$-convergence then yields as $m\to\infty$
\[
\E|u_0c+u_1\xi_1+\cdots+u_n
\xi_n|=\E|u_{i_0}c+u_{i_1}\xi_1+\cdots
+u_{i_n}\xi_n|,
\]
so that $(c,\xi)$ is swap-invariant. Its exchangeability follows
from \cite{har81}, Theorem 1.1.
\end{pf}

\begin{example}[(Example \ref{exlognsi} continuation)]
We show that $n^{-1}(\eta_1+\cdots+\eta_n)$ converges a.s. to
\[
X=\exp\Biggl(\sum_{i=1}^\infty
b_i Z_i - \frac{1}{2} \sum
_{i=1}^\infty b_i^2 \Biggr).
\]
By \cite{kal05s}, Corollary 1.6, we can consider the tail $\sigma$-field
$\sT_{\tilde{\eta}}$, where
\[
\tilde{\eta}=\biggl(\frac{\eta_2}{\eta_1},\frac{\eta_3}{\eta_1},\ldots
\biggr)=
\bigl(\RMe^{Z_2-Z_1-(b_2-b_1)},\RMe^{Z_3-Z_1-(b_3-b_1)},\ldots\bigr).
\]
Since the functions $x \mapsto \RMe^{x-(b_i-b_1)}$, $i\ge2$, are
bijective, $\sT_{\tilde{\eta}}$ can be written as
$\sT_{\tilde{\eta}}=\bigcap_{n\ge2}\mathcal{F}_n$, where
$\mathcal{F}_n=\sigma(Z_n-Z_1,Z_{n+1}-Z_1,\ldots)$. For each $n\ge2$,
the random variable
\[
\tilde{Z}_n=\lim_{k\to\infty}k^{-1} \sum
_{i=0}^{k-1}(Z_1-Z_{n+i})
\]
is clearly $\mathcal{F}_n$-measurable and by the strong law of large
numbers $\tilde{Z}_n=Z_1$ a.s. Thus, $Z_1$ is measurable with
respect to the completion $\bar{\mathcal{F}}_n$ of $\mathcal{F}_n$
for all $n\ge2$, and hence $\bar{\sT}_{\tilde{\eta}}$ measurable. On
the other hand, for all $n\ge2$, the vector $(Z_2,\ldots,Z_n)$ is
independent of $\mathcal{F}_{n+1}$ and therefore independent of
$\sT_{\tilde{\eta}}$. Let $f\dvtx \R\to\R$ be continuous and
bounded. Then for all $A \in\sT_{\tilde{\eta}}$, the dominated
convergence theorem yields
\begin{eqnarray*}
\E\one_A f \Biggl(\sum_{i=2}^\infty
b_i Z_i \Biggr) &=& \lim_{k\to\infty} \E
\one_A f \Biggl(\sum_{i=2}^k
b_iZ_i \Biggr)
\\
&=& \lim_{k\to\infty} \P(A)\E f \Biggl(\sum
_{i=2}^k b_i Z_i \Biggr)=
\P(A)\E f \Biggl(\sum_{i=2}^\infty
b_i Z_i \Biggr),
\end{eqnarray*}
which shows the independence of $\sum_{i=2}^\infty b_i Z_i$ and
$\sT_{\tilde{\eta}}$. Since $\E(Z \mid
\sT_{\tilde{\eta}})=\E(Z\mid\bar{\sT}_{\tilde{\eta}})$ a.s.
for all
integrable $Z$,
\begin{eqnarray*}
\E(\eta_1\mid\sT_{\tilde{\eta}}) &=& \RMe^{(1+b_1)Z_1}
\RMe^{-(1+b_1^2+2b_1)/2},
\\
\E(\eta_2\mid\sT_{\tilde{\eta}}) &=& \RMe^{b_1Z_1}
\RMe^{-b_1^2/2}.
\end{eqnarray*}
By Proposition \ref{proplimit-pos},
\[
\frac{1}{n}\sum_{i=1}^n
\eta_i \to\frac{X \RMe^{Z_1}\RMe^{-(1+2b_1)/2}}{\RMe^{(1+b_1)Z_1}
\RMe^{-(1+b_1^2+2b_1)/2}}\RMe^{b_1Z_1} \RMe^{-b_1^2/2}=X
\qquad\mbox{a.s. and in } L^1 \mbox{ as } n \to\infty.
\]
\end{example}

\section{Non-centred zonoids and lift swap invariance}
\label{secnon-centred-zonoids}

It is possible to relate the centred and non-centred zonoids as
$Z_\xi^o=Z_\xi+Z_{-\xi}$, that is, the centred zonoid is the Minkowski
(elementwise) sum of the zonoid of $\xi$ and the zonoid of $-\xi$
being the central symmetric version $Z_{-\xi}=\{-x\dvt  x\in Z_\xi\}$ of
$Z_\xi$. If $\xi$ has a symmetric distribution, then
$Z_\xi^o=2Z_\xi$ is a scaled zonoid of $\xi$. For a general
integrable $\xi$, its centred zonoid equals $2Z_{\eps\xi}$, where
$\eps$ is the Rademacher random variable taking values $\pm1$ with
equal probability and independent of $\xi$. Note that the conventional
symmetrisation $\xi-\xi'$ for i.i.d. $\xi$ and $\xi'$ is not helpful
in this context.

\begin{proposition}
\label{propn-cent}
If $\xi$ and $\xi^*$ are two integrable random vectors, then
$Z_{\xi}=Z_{\xi^*}$ if and only if $\E\xi=\E\xi^*$ and
$Z_{\xi}^o=Z_{\xi^*}^o$.
\end{proposition}
\begin{pf}
Since $2a_+=|a|+a$ for any real $a$,
\[
h_{Z_\xi}(u) =\tfrac{1}{2}\bigl(\E\bigl|\langle\xi,u\rangle\bigr|+\langle\E
\xi,u\rangle\bigr).
\]
It remains to note that the equality $Z_{\xi}=Z_{\xi^*}$ implies the
equality of expectations by \cite{mos02}, Proposition 2.11.
\end{pf}

In view of the above fact, Proposition \ref{prophf-max} implies that
for positive random vectors the equivalences of centred and
non-centred zonoids are identical concepts.

The centred zonoid of $(1,\xi)$ (also called the centred lift zonoid
of $\xi$) determines uniquely the distribution of $\xi$ by
\cite{har81}, Theorem 1.1. In particular, the invariance of
$\E|1+u_1\xi_1+\cdots+u_d\xi_d|$ with respect to permutations of any
$u_1,\ldots,u_d$ is equivalent to the exchangeability of $\xi$. If the
lifted random vector $(1,\xi)$ is swap-invariant,
that is, $\E|u_0+u_1\xi_1+\cdots+u_d\xi_d|$ is invariant for all
permutations of $u_0,u_1,\ldots,u_d$, then $\xi$ is called \emph{lift
swap-invariant}.

The lift swap-invariance property is slightly weaker than the
\emph{joint self-duality} of $\xi$ meaning the permutation invariance
of $\E(u_0+u_1\xi_1+\cdots+u_d\xi_d)_+$ for all
$u_0,u_1,\ldots,u_d$. The relation between these two properties is
exactly the same as the relation between the equality of centred and
non-centred zonoids. For instance, the lift swap-invariance implies
that $\E|\xi_1|=\cdots=\E|\xi_d|=1$, while the joint self-duality
yields that $\E\xi_1=\cdots=\E\xi_d=1$. The both properties are
identical for random vectors with positive components.

By construction, the lift swap-invariance property implies the
exchangeability of $\xi$ and is actually much stronger. For instance a
vector of i.i.d. positive random variables is exchangeable, but is
neither jointly self-dual nor is lift swap-invariant unless all random
variables equal $1$ almost surely, see \cite{molsch10}.

A weaker version of the self-duality property corresponding to the
permutation of the lifting (constant) coordinate and \emph{one} fixed
other coordinate was studied in \cite{molsch10}. In particular, its
univariate version is often called \emph{put-call symmetry} and is
intensively discussed and applied in the financial literature, see,
for example, \cite{carlee08,teh09} and further references cited
in \cite{molsch10}.

\begin{proposition}
\label{propsd-posit}
If a non-trivial random vector $\xi$ is either jointly self-dual or
is lift swap-invariant with $\E\xi_i=1$ for any $i$, then all its
components are almost surely positive random variables with
expectation being one.
\end{proposition}
\begin{pf}
It suffices to prove this for random variable $\xi$. If $(1,\xi)$ is
swap-invariant and $\E\xi=1$, then $(1,\xi)$ is jointly self-dual by
Proposition \ref{propn-cent}, so it suffices to consider only the
case of a self-dual~$\xi$. The self-duality property of $\xi$
implies that
\[
\E\bigl(0+(-1)\xi\bigr)_+=\E(-1+0\xi)_+,
\]
so that $\E\xi_-=0$ and so $\xi$ is almost surely non-negative.
Since
\[
\E(0+1\xi)_+=\E(1+0\xi)_+,
\]
it follows that $1=\E\xi_+=\E\xi$.

If $\xi$ has an atom at zero, then $\E(1-a\xi)_+$, $a\in\R$, is
bounded from below by a positive number. The self-duality implies
that $\E(-a+\xi)_+$ is also bounded from below by the same number,
which is not possible for large $a$ in view of the integrability of
$\xi$.
\end{pf}

For integrable random vectors with positive components the symmetry
properties can be related to each other. Following the notation of
\cite{molsch11}, define functions
\[
\kappat_j(x)= \biggl(\frac{x_1}{x_j},\ldots,\frac{x_{j-1}}{x_j},
\frac{x_{j+1}}{x_j},\ldots,\frac{x_d}{x_j} \biggr),\qquad j=1,\ldots,d,
\]
on $x\in(0,\infty)^d$. For any $j=1,\ldots,d$ define a new probability
measure by
%
\begin{equation}
\label{eq4} \frac{\mrmdd\P^j}{\mrmdd\P}=\frac{\eta_j}{\E\eta_j}.
\end{equation}
This measure change was used in \cite{eberpapshir08b} in order to
reduce the dimensionality when calculating option prices. Consider an
integrable random vector $\eta$ with positive components. If
$\E\eta_j=1$, then the zonoid of $\eta$ coincides with the lift zonoid
of $\kappat_j(\eta)$ under $\P^j$, see \cite{molsch11m},
Proposition 3.

\begin{theorem}
\label{proprelations}
Assume that $\eta$ is an integrable random vector of dimension
$d\geq2$ with positive components. The following conditions are
equivalent:
\begin{itemize}[(c)]
\item[(a)] $\eta$ is swap-invariant under $\P$.
\item[(b)] $\kappat_j(\eta)$ is lift swap-invariant (equivalently
jointly self-dual) under $\P^j$ for any (and then all)
$j\in\{1,\ldots,d\}$.
\item[(c)] In case $d\geq3$, for at least two
$j\in\{1,\ldots,d\}$ (and then automatically for all $j$),
$\kappat_j(\eta)$ is exchangeable under $\P^j$.
\end{itemize}
\end{theorem}
\begin{pf}
The equivalence of (a) and (b) is obtained (for $j=1$)
by
\[
\E|u_1\eta_1+\cdots+u_d\eta_d| =
\E\eta_1 \E_{\P^1} \biggl|u_1+u_2
\frac{\eta_2}{\eta_1}+\cdots+u_d \frac{\eta_d}{\eta_1} \biggr|,
\]
so that permutations of coordinates in the left-hand side
corresponds to permutations in the right-hand side. The invariance
with respect to the latter is equivalent to the lift swap invariance
of $\tilde{\kappa}_1(\eta)$ under $\P^1$, since the right-hand side
identifies the distribution of $\tilde{\kappa}_1(\eta)$.

It is easy to see that (a) implies (c) for all $j$,
since the exchangeability is a weaker property than
(b). Assuming (c) for $j=1,2$ without loss of
generality, we see that $(\eta_2/\eta_1,\ldots,\eta_d/\eta_1)$ is
$\P^1$-exchangeable and
$(\eta_1/\eta_2,\eta_3/\eta_2,\ldots,\eta_d/\eta_2)$ is
$\P^2$-exchangeable. The first fact implies that $\E|\langle
u,\eta\rangle|$ is invariant with respect to permutation all but
first coordinates of $u$, while the second fact implies the
invariance with respect to permutations of all coordinates excluding
the second one, so $\eta$ is swap-invariant.
\end{pf}

\section{Equality of zonoids}
\label{secequality-zonoids}

\subsection{Location-scale families}
\label{seclocat-scale-famil}

Consider family of random variables $\xi=\mu+\sigma X$ for an
integrable random variable $X$ and $\mu\in\R$, $\sigma>0$. These random
variables are said to form a location-scale family.

\begin{theorem}
\label{thrloc-scale}
Assume that the distribution of $X$ has infinite essential infimum
and essential supremum. Then the zonoid $Z_\xi$ of a random variable
$\xi$ from the location-scale family generated by $X$ uniquely
determines the location and scale parameters of the distribution.
\end{theorem}
\begin{pf}
Without loss of generality, set $\E X=0$. Assume that the random
variables $\mu+\sigma X$ and $\mu^*+\sigma^*X$ share the same
zonoid. By Proposition \ref{propn-cent}, $\mu=\mu^*$.

In order to finish the proof, we show that $\E(\mu+\sigma X)_+$ is
strictly increasing in $\sigma$ for each fixed $\mu\in\R$. This is
obvious if $\mu=0$, since $\E(\sigma X)_+=\sigma\E X_+$, which is
strictly increasing in $\sigma$ since $\E X_+>0$.

Assume that $\mu<0$ and $\sigma_1 > \sigma_2$. Then
\begin{eqnarray*}
&&\E\bigl((\mu+\sigma_1X)_+-(\mu+\sigma_2X)_+\bigr)
\\
&&\quad=\E\bigl((\mu+\sigma_1X)\one_{\{-{\mu}/{\sigma_1} < X \le
-{\mu}/{\sigma_2}\}}\bigr) +(
\sigma_1-\sigma_2)\E(X\one_{\{-{\mu}/{\sigma_2}< X\}})>0,
\end{eqnarray*}
where the last expectation is positive because $X$
has unbounded support and $\E X=0$.

If $\mu> 0$, the same argument applied to
$\E(\mu+\sigma_1X)_-$ yields that the expectation of the negative
part is strictly decreasing in $\sigma$ and the equality
$\E(\mu+\sigma_1X)_+=\mu-\E(\mu+\sigma_1X)_-$ concludes the proof.
\end{pf}

Note that Theorem \ref{thrloc-scale} does not hold for the centred
zonoid $Z^o_\xi$ unless it is assumed that the expectation of $\xi$ is
known and so $Z_\xi$ is also identified.

\begin{corollary}
\label{corlin}
Assume that random variable $\xi$ has infinite essential infimum and
essential supremum. If $Z_\xi=Z_{\sigma\xi+\mu}$, then $\mu=0$ and
$\sigma=1$.
\end{corollary}

\begin{corollary}
\label{coruniquely-normal}
Two normally distributed $d$-dimensional random vectors $\xi$ and $\xi
^*$ coincide
in distribution if and only if $Z_\xi=Z_{\xi^*}$.
\end{corollary}
\begin{pf}
For $u\in\R^d$ the random variables $\langle\xi,u\rangle$ and
$\langle\xi^*,u\rangle$ belong to the same location-scale
family. The proof is finished by referring to
Theorem \ref{thrloc-scale} and noticing that all one-dimensional
projection of a random vector uniquely determine its distribution.
\end{pf}

The uniqueness holds also for the location scale family obtained as
$\mu+\sigma X$ for a symmetric stable random variable $X$.

\begin{example}[(Distribution with bounded support)]
Assume that $\E X=0$ and that $X$ has finite essential infinum,
that is, there exists a constant $c$ such that $X\ge c$ a.s. Choose
$\mu
>0$. Then for all $\sigma< -\mu/c$ the random variable
$\xi=\mu+\sigma X$ is a.s. positive and so the expectation of its
negative part is zero and the expectation of its positive part is
$\mu$. Thus, the zonoid $Z_\xi$ does not uniquely determine the scale
parameter $\sigma$.
\end{example}

Note that all above results are formulated for non-centred zonoids. In
the rest of this section, we consider centred zonoids, and the
corresponding zonoid equivalence concept. The following result
concerns random vectors that can be represented as product of a
scaling random variable and an independent random vector.

\begin{proposition}
\label{propscaling}
Two random vectors $\xi=R\zeta$ and $\xi^*=R^*\zeta^*$, where $R$
and $R^*$ are positive random variables independent of $\zeta$
and $\zeta^*$, respectively, are zonoid equivalent if and only if
$(\E R)\zeta$ and $(\E R^*)\zeta^*$ are zonoid equivalent.
\end{proposition}
\begin{pf}
It suffices to note that
\[
\E\bigl|\langle u,\xi\rangle\bigr|=\E R\E\bigl|\langle u,\zeta\rangle\bigr| =\E\bigl|\bigl
\langle u,(\E
R)\zeta\bigr\rangle\bigr|.
\]
\upqed
\end{pf}

Consider random vectors with centred \emph{elliptical} distributions,
that is, assume that $\xi=R(AU)$, where $U$ is uniformly distributed on
the unit sphere, $A$ is a (deterministic) matrix and $R$ is a positive
random variable independent of $U$.

\begin{proposition}
\label{propellip}
Two centred elliptically distributed random vectors $\xi=R(AU)$ and
$\xi^*=R^*(A^*U)$ are zonoid equivalent if and only if $(\E R)^2
AA^\top=(\E R^*)^2 A^*(A^*)^\top$.
\end{proposition}
\begin{pf}
Using rescaling, it is possible to assume that $\E R=\E R^*$. By
Proposition \ref{propscaling}, it suffices to consider zonoid
equivalence of $AU$ and $A^*U$. By Proposition \ref{proplinear},
this is the case if and only if random variables $\langle A^\top
u,U\rangle$ and $\langle(A^*)^\top u,U\rangle$ are zonoid
equivalent. Since $U$ is uniformly distributed on the unit sphere,
$\langle v,U\rangle$ is distributed as a certain random variable
with a fixed distribution scaled by $\|v\|$ for all $v$. Thus,
$\|A^\top u\|=\|(A^*)^\top u\|$ for all $u$, which implies
the statement.
\end{pf}

\begin{corollary}
\label{cormultivariate-normal}
Two symmetric normally distributed random vectors $\xi$ and $\xi^*$
coincide in distribution if and only if they are zonoid equivalent.
\end{corollary}

Zonoid of $S\alpha S$ random $\xi$ with $\alpha\in(1,2]$ is computed
in \cite{mo09}, Section 6.4, as
\[
Z_\xi=\frac{1}{\uppi}\Gamma\biggl(1-\frac{1}{\alpha}\biggr)K,
\]
where $\Gamma$ is the gamma-function and $K$ is a convex body that,
together with $\alpha$, characterises the distribution of $\xi$. Thus,
if $\alpha$ is fixed, then the zonoid determines uniquely the
corresponding symmetric $\alpha$-stable distribution. However, two
symmetric stable vectors with the same zonoid are not necessarily
identically distributed if their stability indices are different.

\subsection{Log-infinitely divisible distributions with equal zonoids}
\label{seclog-inf}

A random vector with positive components can be written as the
coordinate-wise exponential $\eta=\RMe^{\xi}$. In the following,
$\varphi_\xi$ stands for the characteristic function of $\xi$. The
following result immediately follows from
\cite{kabschhaan09}, Proposition 6, see also \cite{molsch11}, Theorem 3.2.

\begin{theorem}
\label{thchar-equal-zonoids}
Two integrable random vectors $\RMe^{\xi}$ and $\RMe^{\xi^*}$ are zonoid
equivalent if and only if
%
\begin{equation}
\label{eqchar-f-zon} \varphi_\xi(u-\imagi w)=\varphi_{\xi^*}(u-
\imagi w)
\end{equation}
for all $u\in\R^d$ with $\sum u_i=0$ and for at least one (and then
necessarily for all) $w$, such that $\sum w_k=1$ and both sides in
(\ref{eqchar-f-zon}) are finite.
\end{theorem}

Assume that $\RMe^\xi$ and $\RMe^{\xi^*}$ are two random vectors, where
$\xi$ and $\xi^*$ are infinitely divisible random variables. Then
\[
\varphi_\xi(u)=\E \RMe^{\imagi\langle u,\xi\rangle}= \exp\biggl\{\imagi
\langle b,u
\rangle-\frac{1}{2}\langle u,Au\rangle+\int_{\R^d}
\bigl(\RMe^{\imagi\langle u,x\rangle}-1 -\imagi\langle u,x\rangle\one_{\|
x\|\leq1}\bigr)\mrmd
\nu(x) \biggr\}
\]
for $u\in\R^d$, where $A=(a_{ij})$ is a symmetric non-negative definite
$d\times d$ matrix, $b\in\R^d$ is a constant vector and $\nu$ is a
measure on $\R^d$ (called the L\'evy measure) satisfying
$\nu(\{0\})=0$ and
\[
\int_{\R^d}\min\bigl(\|x\|^2,1\bigr)\mrmd \nu(x)<
\infty.
\]
Then $\xi$ is said to have the L\'evy triplet $(A,\nu,b)$.
In this section, we translate the equality of the zonoids of two
log-infinitely divisible random vectors into conditions on their
L\'evy triplets. Note that the conditions on the L\'evy triplet of
infinitely divisible random vectors apply also for L\'evy processes
with time one values $\xi$ and $\xi^*$.

In order to formulate the condition on the Gaussian terms in a compact
form it is helpful to use the \emph{variogram}
\[
\gamma_{ij}=a_{ii}+a_{jj}-2a_{ij}.
\]
If $\xi$ is normally distributed, then $\gamma_{ij}$ is the variance
of $\xi_i-\xi_j$. In order to state the condition on the L\'evy
measure define $(d-1)\times d$-dimensional matrix, $d\ge2$
%
\begin{equation}
\label{equ} U=\pmatrix{ 1 & 0 & \cdots& 0 & -1
\cr
0 & 1 & \cdots& 0 & -1
\cr
\vdots& \vdots& \ddots& \vdots& \vdots
\cr
0 & 0 & \cdots& 1 & -1}.
\end{equation}

\begin{theorem}
\label{thchar-equal-zonoids-inf-div}
Let $\RMe^{\xi}$ and $\RMe^{\xi^*}$ be integrable random
vectors such that $\xi$ and $\xi^*$ are infinitely divisible with
characteristic triplets $(A,\nu,\gamma)$ and $(A^*,\nu^*,\gamma^*)$.
Then for $d\geq2$ $\RMe^{\xi}$ and $\RMe^{\xi^*}$ are zonoid equivalent if
and only if the
following three conditions hold.
\begin{itemize}[(c)]
\item[(a)] $\gamma_{ij}=\gamma^*_{ij}$ for all
$i,j\in\{1,\ldots,d\}$.
\item[(b)] The images $\hat\nu U^{-1}$ and $\hat\nu^* U^{-1}$ under
$U$ of measures $\mrmdd\hat\nu(x)=\RMe^{x_d}\mrmd \nu(x)$ and
$\mrmdd\hat\nu^*(x)=\RMe^{x_d}\mrmd \nu^*(x)$, $x\in\R^d$, restricted to
$\R^{d-1}\setminus\{0\}$ coincide.
\item[(c)] $\E \RMe^{\xi_i}=\E \RMe^{\xi^*_i}$ for all $i=1,\ldots,d$, that
is,
%
\begin{eqnarray}
\label{eqbexp}
&&
b_i+\frac{1}{2} a_{ii}+\int
_{\R^d}\bigl(\RMe^{x_i}-1-x_i
\one_{\|x\|\leq1}\bigr)\mrmd \nu(x)
\nonumber\\[-8pt]\\[-8pt]
&&\quad=b_i^*+\frac{1}{2} a_{ii}^*+\int
_{\R^d} \bigl(\RMe^{x_i}-1-x_i
\one_{\|x\|\leq1}\bigr)\mrmd \nu^*(x).
\nonumber
\end{eqnarray}
\end{itemize}
For $d=1$, $\RMe^{\xi}$ and $\RMe^{\xi^*}$ are zonoid equivalent if and only
if \textup{(c)} holds.
\end{theorem}

The following result is closely related to and can be alternatively
derived following the proof of \cite{kabschhaan09}, Theorem 10, see also
\cite{kab10}, Theorem 1.1.

\begin{corollary}
\label{corlognorm}
Two lognormal random vectors $\RMe^\xi$ and $\RMe^{\xi^*}$ are zonoid
equivalent if and only if $\mu_i+\frac{1}{2} a_{ii}=\mu^*_i+\frac{1}{2}
a^*_{ii}$
for all $i$ and $\gamma_{ij}=\gamma_{ij}^*$ for all $i,j$,
that is, $\xi$ and $\xi^*$ have identical variogram.
\end{corollary}

In particular, in the lognormal case the zonoid equivalence does not
even imply the equality of the marginal distributions, quite
differently to the case of normal distributions where the zonoid
uniquely determines the joint distribution, see
Corollary \ref{cormultivariate-normal}.

Furthermore, note that the kernel of $U$ given by (\ref{equ}) is the
family of vectors with all equal components. Hence, if the support of
$\nu$ is a subset of the kernel of $U$, then the corresponding
log-infinitely divisible distribution shares the same zonoid with a
lognormal distribution, meaning that two rather different
distributions are zonoid equivalent.

\begin{pf*}{Proof of Theorem \ref{thchar-equal-zonoids-inf-div}}
For $d\geq2$, the zonoid equivalence of $\RMe^\xi$ and $\RMe ^{\xi^*}$
implies $\E \RMe ^\xi=\E \RMe ^{\xi^*}$, see Proposition \ref{prophf-max},
and in particular $c=\E \RMe ^{\xi_d}=\E \RMe ^{\xi^*_d}$. Note that this is
also implied by (c). Since also $Z_{\RMe ^\xi}=Z_{\RMe ^{\xi^*}}$ by
Proposition \ref{prophf-max},
\[
\E\bigl(u_1\RMe ^{\xi_1}+\cdots+u_d\RMe^{\xi_d}
\bigr)_+ =\E \RMe ^{\xi_d}\bigl(u_1\RMe ^{\xi_1-\xi_d}+\cdots+u_{d-1}\RMe ^{\xi_{d-1}-\xi
_1}+u_d\bigr)_+,
\]
the zonoid of $\RMe ^\xi$ uniquely determines and is uniquely determined
by the probability distribution of
$U\xi=(\xi_1-\xi_d,\ldots,\xi_{d-1}-\xi_d)$ under the probability
measure $\P^d$ with density $\RMe ^{\xi_d}/c$.

In order to identify the distribution of $U\xi$ under $\P^d$ first
note that the distribution of $\xi$ under $\P^d$ has the
characteristic triplet $(A,\hat\nu,\hat{b})$, where
$\mrmdd\hat\nu(x)=\RMe ^{x_d}\mrmd \nu(x)$ and
\[
\hat b=b+\int_{\|x\|\leq1}x\bigl(\RMe ^{x_d}-1\bigr)
\nu(\mrmdd x)+Ae_{d},
\]
where $e_d$ is the $d$th standard basis vector, see \cite{sat00}, Example 7.3. By \cite{sat99}, Proposition 11.10, the L\'evy
triplet of $U\xi$ under $\P^d$ is given by $A_U=UAU^\top$, $\hat\nu
U^{-1}$ restricted onto $\R^{d-1}\setminus\{0\}$ and
\[
b_U=U\hat b+\int_{\R^d}Ux(\one_{\|Ux\|\leq
1}-
\one_{\|x\|\leq1})\hat\nu(\mrmdd x).
\]
The corresponding formula holds for $\xi^*$.

Equating the centred Gaussian terms, the L\'evy measures, and
simplifying $b_U=b_U^*$ yields that $U\xi$ under $\P^d$ coincides in
distribution with $U\xi^*$ under $\P^{d*}$ if and only if
%
\begin{equation}
\label{eqasi-two}
a_{ij}+a_{dd}-a_{di}-a_{jd}=a_{ij}^*+a_{dd}^*-a_{di}^*-a_{jd}^*,\qquad
i,j=1,\ldots,d-1,
\end{equation}
condition (b) holds and, for all $i=1,\ldots,d-1$,
%
\begin{eqnarray}
\label{eqbb}
&&
b_i-b_d+a_{id}-a_{dd}+
\int_{\R^d}(x_i-x_d) \bigl(
\one_{\|Ux\|\leq
1}\RMe ^{x_d}-\one_{\|x\|\leq1}\bigr)\mrmd \nu(x)
\nonumber\\[-8pt]\\[-8pt]
&&\quad=b_i^*-b_d^*+a_{id}^*-a_{dd}^*+
\int_{\R^d}(x_i-x_d) \bigl(
\one_{\|Ux\|\leq
1}\RMe ^{x_d}-\one_{\|x\|\leq1}\bigr)\mrmd \nu^*(x).
\nonumber
\end{eqnarray}

Adding equations (\ref{eqasi-two}) with $k,l=i,i$; $k,l=j,j$ (for
given $i$ and $j$), and subtracting (\ref{eqasi-two}) multiplied by
two, we arrive at the equality of the variograms. Furthermore,
noticing that
\[
(a_{ij}+a_{dd}-a_{di}-a_{jd})_{ij=1}^{d-1}
=\tfrac{1}{2}(\gamma_{id}+\gamma_{jd}-
\gamma_{ij})_{ij=1}^{d-1}
\]
we obtain that the equality of variograms implies
(\ref{eqasi-two}). The equality of zonoids implies the equality of
expectations, which exactly corresponds to (\ref{eqbexp}). It
remains to show that (\ref{eqbexp}) together with other two
conditions (a) and (b) imply (\ref{eqbb}).

By (\ref{eqbexp}), we have for all $i=1,\ldots,d-1$
%
\begin{eqnarray}
\label{eqequal-exp-i}
&&
b_i+\frac{1}{2} a_{ii}+\int_{\R^d}
\bigl(\RMe ^{x_i}-1-x_i\one_{\|x\|\leq
1}\bigr)\mrmd \nu(x)
\nonumber\\[-8pt]\\[-8pt]
&&\quad=b_i^*+\frac{1}{2}
a_{ii}^*+\int_{\R^d}\bigl(\RMe ^{x_i}-1-x_i
\one_{\|x\|\leq
1}\bigr)\mrmd \nu^*(x),
\nonumber\\
\label{eqequal-exp-n}
&&
b_d+\frac{1}{2} a_{dd}+\int_{\R^d}
\bigl(\RMe ^{x_d}-1-x_d\one_{\|x\|\leq
1}\bigr)\mrmd \nu(x)
\nonumber\\[-8pt]\\[-8pt]
&&\quad=b_d^*+\frac{1}{2} a_{dd}^*+
\int_{\R^d}\bigl(\RMe ^{x_d}-1-x_d
\one_{\|x\|\leq
1}\bigr)\mrmd \nu^*(x),\nonumber
\end{eqnarray}
while condition (a) implies
%
\begin{equation}
\label{eqmat-ii} a_{ii}+a_{dd}-2a_{id}=a_{ii}^*+a_{dd}^*-2a_{id}^*
\end{equation}
for all $i=1,\ldots,d-1$. Furthermore, condition (b) implies
%
\begin{eqnarray}
\label{eqeq-integrals}
&&
\int_{\R^d}\bigl(\RMe ^{x_i-x_d}-1-(x_i-x_d)
\one_{\|Ux\|\leq1}\bigr)\mrmd \hat\nu(x)
\\
&&\quad=\int_{\R^d}\bigl(\RMe ^{x_i-x_d}-1-(x_i-x_d)
\one_{\|Ux\|\leq1}\bigr)\mrmd \hat\nu^*(x),
\nonumber
\end{eqnarray}
where $d\hat\nu(x)=\RMe ^{x_d}\mrmd \nu(x)$, since by changing variables
\[
\int_{\R^{d-1}}\bigl(\RMe ^y-1-y\one_{\|y\|\leq1}
\bigr)\mrmd \bigl(\hat\nu U^{-1}\bigr) (y) =\int_{\R^{d-1}}
\bigl(\RMe ^y-1-y\one_{\|y\|\leq1}\bigr)\mrmd \bigl(\hat\nu^*
U^{-1}\bigr) (y).
\]
Now (\ref{eqbb}) is obtained by subtracting
from (\ref{eqequal-exp-i}) the sum
of (\ref{eqeq-integrals}), (\ref{eqequal-exp-n}) and a half
of (\ref{eqmat-ii}).

Recall that equality of the zonoids is equivalent to equality of
their support functions for all $u$ on the unite sphere. Hence, for
positive random variables $\RMe ^\xi$ and $\RMe ^{\xi^*}$ ($d=1$) equality
of their zonoids is equivalent to equality of their expectations,
which in turn, is equivalent to condition~(c).
\end{pf*}

\section*{Acknowledgements}

This work was supported by Swiss National Science Foundation Grants
200021-126503 and 200021-137527 and has been finished while IM held
the Chair of Excellence at the University Carlos III of Madrid
supported by the Santander bank.

The authors are grateful to Markus Kiderlen for useful information
concerning zonoids. The thoughtful comments and constructive
suggestions of the Associate Editor and the referees have led to a
clarification of the exposition and a proper accentuation of
relationships with the theory of stable laws.


%

\printhistory


\begin{thebibliography}{41}

\bibitem{bor09}
\begin{barticle}[mr]
\bauthor{\bsnm{Borell},~\bfnm{Christer}\binits{C.}}
(\byear{2009}).
\btitle{Zonoids induced by {G}auss measure with an application to risk
  aversion}.
\bjournal{ALEA Lat. Am. J. Probab. Math. Stat.}
\bvolume{6}
\bpages{133--147}.
\bid{issn={1980-0436}, mr={2506861}}
\bptok{imsref}%
\end{barticle}
\endbibitem

\bibitem{brelit78}
\begin{barticle}[author]
\bauthor{\bsnm{Breeden},~\bfnm{D.~T.}\binits{D.T.}} \AND
  \bauthor{\bsnm{Litzenberger},~\bfnm{R.~H.}\binits{R.H.}}
(\byear{1978}).
\btitle{Prices of state-contingent claims implicit in options prices}.
\bjournal{J. Business}
\bvolume{51}
\bpages{621--651}.
\bptok{imsref}%
\end{barticle}
\endbibitem

\bibitem{brores77}
\begin{barticle}[mr]
\bauthor{\bsnm{Brown},~\bfnm{Bruce~M.}\binits{B.M.}} \AND
  \bauthor{\bsnm{Resnick},~\bfnm{Sidney~I.}\binits{S.I.}}
(\byear{1977}).
\btitle{Extreme values of independent stochastic processes}.
\bjournal{J. Appl. Probab.}
\bvolume{14}
\bpages{732--739}.
\bid{issn={0021-9002}, mr={0517438}}
\bptok{imsref}%
\end{barticle}
\endbibitem

\bibitem{carellgup98}
\begin{barticle}[author]
\bauthor{\bsnm{Carr},~\bfnm{P.}\binits{P.}},
  \bauthor{\bsnm{Ellis},~\bfnm{K.}\binits{K.}} \AND
  \bauthor{\bsnm{Gupta},~\bfnm{V.}\binits{V.}}
(\byear{1998}).
\btitle{Static hedging of exotic options}.
\bjournal{J. Finance}
\bvolume{53}
\bpages{1165--1190}.
\bptok{imsref}%
\end{barticle}
\endbibitem

\bibitem{carlee08}
\begin{barticle}[mr]
\bauthor{\bsnm{Carr},~\bfnm{Peter}\binits{P.}} \AND
  \bauthor{\bsnm{Lee},~\bfnm{Roger}\binits{R.}}
(\byear{2009}).
\btitle{Put-call symmetry: Extensions and applications}.
\bjournal{Math. Finance}
\bvolume{19}
\bpages{523--560}.
\bid{doi={10.1111/j.1467-9965.2009.00379.x}, issn={0960-1627}, mr={2583519}}
\bptok{imsref}%
\end{barticle}
\endbibitem

\bibitem{cas10}
\begin{bincollection}[mr]
\bauthor{\bsnm{Cascos},~\bfnm{Ignacio}\binits{I.}}
(\byear{2010}).
\btitle{Data depth: Multivariate statistics and geometry}.
In \bbooktitle{New Perspectives in Stochastic Geometry}
(\beditor{\bfnm{W.~S.}\binits{W.S.}~\bsnm{Kendall}} \AND
  \beditor{\bfnm{I.}\binits{I.}~\bsnm{Molchanov}}, eds.)
\bpages{398--423}.
\blocation{Oxford}: \bpublisher{Oxford Univ. Press}.
\bid{mr={2654685}}
\bptok{imsref}%
\end{bincollection}
\endbibitem

\bibitem{dac82}
\begin{bincollection}[mr]
\bauthor{\bsnm{Dacunha-Castelle},~\bfnm{D.}\binits{D.}}
(\byear{1982}).
\btitle{A survey on exchangeable random variables in normed spaces}.
In \bbooktitle{Exchangeability in Probability and Statistics ({R}ome, 1981)}
(\beditor{\bfnm{G.}\binits{G.}~\bsnm{Koch}}
\AND \beditor{\bfnm{F.}\binits{F.}~\bsnm{Spizzichino}}, eds.)
\bpages{47--60}.
\blocation{Amsterdam}: \bpublisher{North-Holland}.
\bid{mr={0675964}}
\bptok{imsref}%
\end{bincollection}
\endbibitem

\bibitem{davj}
\begin{bbook}[mr]
\bauthor{\bsnm{Daley},~\bfnm{D.~J.}\binits{D.J.}} \AND
  \bauthor{\bsnm{Vere-Jones},~\bfnm{D.}\binits{D.}}
(\byear{1988}).
\btitle{An Introduction to the Theory of Point Processes}.
\bseries{Springer Series in Statistics}.
\blocation{New York}: \bpublisher{Springer}.
\bid{mr={0950166}}
\bptok{imsref}%
\end{bbook}
\endbibitem

\bibitem{haa84}
\begin{barticle}[mr]
\bauthor{\bparticle{de} \bsnm{Haan},~\bfnm{L.}\binits{L.}}
(\byear{1984}).
\btitle{A spectral representation for max-stable processes}.
\bjournal{Ann. Probab.}
\bvolume{12}
\bpages{1194--1204}.
\bid{issn={0091-1798}, mr={0757776}}
\bptok{imsref}%
\end{barticle}
\endbibitem

\bibitem{haapic86}
\begin{barticle}[mr]
\bauthor{\bparticle{de} \bsnm{Haan},~\bfnm{L.}\binits{L.}} \AND
  \bauthor{\bsnm{Pickands},~\bfnm{J.}\binits{J.} \bsuffix{III}}
(\byear{1986}).
\btitle{Stationary min-stable stochastic processes}.
\bjournal{Probab. Theory Related Fields}
\bvolume{72}
\bpages{477--492}.
\bid{doi={10.1007/BF00344716}, issn={0178-8051}, mr={0847381}}
\bptok{imsref}%
\end{barticle}
\endbibitem

\bibitem{eberpapshir08b}
\begin{barticle}[mr]
\bauthor{\bsnm{Eberlein},~\bfnm{Ernst}\binits{E.}},
  \bauthor{\bsnm{Papapantoleon},~\bfnm{Antonis}\binits{A.}} \AND
  \bauthor{\bsnm{Shiryaev},~\bfnm{Albert~N.}\binits{A.N.}}
(\byear{2009}).
\btitle{Esscher transform and the duality principle for multidimensional
  semimartingales}.
\bjournal{Ann. Appl. Probab.}
\bvolume{19}
\bpages{1944--1971}.
\bid{doi={10.1214/09-AAP600}, issn={1050-5164}, mr={2569813}}
\bptok{imsref}%
\end{barticle}
\endbibitem

\bibitem{falhusrei04}
\begin{bbook}[mr]
\bauthor{\bsnm{Falk},~\bfnm{Michael}\binits{M.}},
  \bauthor{\bsnm{H{\"u}sler},~\bfnm{J{\"u}rg}\binits{J.}} \AND
  \bauthor{\bsnm{Reiss},~\bfnm{Rolf-Dieter}\binits{R.D.}}
(\byear{2004}).
\btitle{Laws of Small Numbers: Extremes and Rare Events},
\bedition{extended} ed.
\blocation{Basel}: \bpublisher{Birkh\"auser}.
\bid{mr={2104478}}
\bptok{imsref}%
\end{bbook}
\endbibitem

\bibitem{har81}
\begin{barticle}[mr]
\bauthor{\bsnm{Hardin},~\bfnm{Clyde~D.}\binits{C.D.} \bsuffix{Jr.}}
(\byear{1981}).
\btitle{Isometries on subspaces of {$L\sp{p}$}}.
\bjournal{Indiana Univ. Math. J.}
\bvolume{30}
\bpages{449--465}.
\bid{doi={10.1512/iumj.1981.30.30036}, issn={0022-2518}, mr={0611233}}
\bptok{imsref}%
\end{barticle}
\endbibitem

\bibitem{har82}
\begin{barticle}[mr]
\bauthor{\bsnm{Hardin},~\bfnm{Clyde~D.}\binits{C.D.} \bsuffix{Jr.}}
(\byear{1982}).
\btitle{On the spectral representation of symmetric stable processes}.
\bjournal{J.~Multivariate Anal.}
\bvolume{12}
\bpages{385--401}.
\bid{doi={10.1016/0047-259X(82)90073-2}, issn={0047-259X}, mr={0666013}}
\bptok{imsref}%
\end{barticle}
\endbibitem

\bibitem{kab09s}
\begin{barticle}[mr]
\bauthor{\bsnm{Kabluchko},~\bfnm{Zakhar}\binits{Z.}}
(\byear{2009}).
\btitle{Spectral representations of sum- and max-stable processes}.
\bjournal{Extremes}
\bvolume{12}
\bpages{401--424}.
\bid{doi={10.1007/s10687-009-0083-9}, issn={1386-1999}, mr={2562988}}
\bptok{imsref}%
\end{barticle}
\endbibitem

\bibitem{kab10}
\begin{barticle}[mr]
\bauthor{\bsnm{Kabluchko},~\bfnm{Zakhar}\binits{Z.}}
(\byear{2010}).
\btitle{Stationary systems of {G}aussian processes}.
\bjournal{Ann. Appl. Probab.}
\bvolume{20}
\bpages{2295--2317}.
\bid{doi={10.1214/10-AAP686}, issn={1050-5164}, mr={2759735}}
\bptok{imsref}%
\end{barticle}
\endbibitem

\bibitem{kabschhaan09}
\begin{barticle}[mr]
\bauthor{\bsnm{Kabluchko},~\bfnm{Zakhar}\binits{Z.}},
  \bauthor{\bsnm{Schlather},~\bfnm{Martin}\binits{M.}} \AND
  \bauthor{\bparticle{de} \bsnm{Haan},~\bfnm{Laurens}\binits{L.}}
(\byear{2009}).
\btitle{Stationary max-stable fields associated to negative definite
  functions}.
\bjournal{Ann. Probab.}
\bvolume{37}
\bpages{2042--2065}.
\bid{doi={10.1214/09-AOP455}, issn={0091-1798}, mr={2561440}}
\bptok{imsref}%
\end{barticle}
\endbibitem

\bibitem{kabstoev12}
\begin{bmisc}[author]
\bauthor{\bsnm{Kabluchko},~\bfnm{Z.}\binits{Z.}} \AND
  \bauthor{\bsnm{Stoev},~\bfnm{S.}\binits{S.}}
(\byear{2012}).
\bhowpublished{Minimal spectral representations of infinitely divisible and
  max-infinitely divisible processes. Technical report.
Available at \arxivurl{arXiv:1207.4983}}.
\bptok{imsref}%
\end{bmisc}
\endbibitem

\bibitem{kalle}
\begin{bbook}[mr]
\bauthor{\bsnm{Kallenberg},~\bfnm{Olav}\binits{O.}}
(\byear{2002}).
\btitle{Foundations of Modern Probability},
\bedition{2nd} ed.
\bseries{Probability and Its Applications (New York)}.
\blocation{New York}: \bpublisher{Springer}.
\bid{mr={1876169}}
\bptok{imsref}%
\end{bbook}
\endbibitem

\bibitem{kal05s}
\begin{bbook}[mr]
\bauthor{\bsnm{Kallenberg},~\bfnm{Olav}\binits{O.}}
(\byear{2005}).
\btitle{Probabilistic Symmetries and Invariance Principles}.
\bseries{Probability and Its Applications (New York)}.
\blocation{New York}: \bpublisher{Springer}.
\bid{mr={2161313}}
\bptok{imsref}%
\end{bbook}
\endbibitem

\bibitem{kosmos98}
\begin{barticle}[mr]
\bauthor{\bsnm{Koshevoy},~\bfnm{Gleb}\binits{G.}} \AND
  \bauthor{\bsnm{Mosler},~\bfnm{Karl}\binits{K.}}
(\byear{1998}).
\btitle{Lift zonoids, random convex hulls and the variability of random
  vectors}.
\bjournal{Bernoulli}
\bvolume{4}
\bpages{377--399}.
\bid{doi={10.2307/3318721}, issn={1350-7265}, mr={1653276}}
\bptok{imsref}%
\end{barticle}
\endbibitem

\bibitem{lepwoodzin81}
\begin{barticle}[mr]
\bauthor{\bsnm{LePage},~\bfnm{Raoul}\binits{R.}},
  \bauthor{\bsnm{Woodroofe},~\bfnm{Michael}\binits{M.}} \AND
  \bauthor{\bsnm{Zinn},~\bfnm{Joel}\binits{J.}}
(\byear{1981}).
\btitle{Convergence to a stable distribution via order statistics}.
\bjournal{Ann. Probab.}
\bvolume{9}
\bpages{624--632}.
\bid{issn={0091-1798}, mr={0624688}}
\bptok{imsref}%
\end{barticle}
\endbibitem

\bibitem{mo1}
\begin{bbook}[mr]
\bauthor{\bsnm{Molchanov},~\bfnm{Ilya}\binits{I.}}
(\byear{2005}).
\btitle{Theory of Random Sets}.
\bseries{Probability and Its Applications (New York)}.
\blocation{London}: \bpublisher{Springer}.
\bid{mr={2132405}}
\bptok{imsref}%
\end{bbook}
\endbibitem

\bibitem{mo08e}
\begin{barticle}[mr]
\bauthor{\bsnm{Molchanov},~\bfnm{Ilya}\binits{I.}}
(\byear{2008}).
\btitle{Convex geometry of max-stable distributions}.
\bjournal{Extremes}
\bvolume{11}
\bpages{235--259}.
\bid{doi={10.1007/s10687-008-0055-5}, issn={1386-1999}, mr={2429906}}
\bptok{imsref}%
\end{barticle}
\endbibitem

\bibitem{mo09}
\begin{barticle}[mr]
\bauthor{\bsnm{Molchanov},~\bfnm{Ilya}\binits{I.}}
(\byear{2009}).
\btitle{Convex and star-shaped sets associated with multivariate stable
  distributions. {I}. {M}oments and densities}.
\bjournal{J. Multivariate Anal.}
\bvolume{100}
\bpages{2195--2213}.
\bid{doi={10.1016/j.jmva.2009.04.003}, issn={0047-259X}, mr={2560363}}
\bptok{imsref}%
\end{barticle}
\endbibitem

\bibitem{molsch10}
\begin{barticle}[mr]
\bauthor{\bsnm{Molchanov},~\bfnm{Ilya}\binits{I.}} \AND
  \bauthor{\bsnm{Schmutz},~\bfnm{Michael}\binits{M.}}
(\byear{2010}).
\btitle{Multivariate extension of put-call symmetry}.
\bjournal{SIAM J. Financial Math.}
\bvolume{1}
\bpages{396--426}.
\bid{doi={10.1137/090754194}, issn={1945-497X}, mr={2652071}}
\bptok{imsref}%
\end{barticle}
\endbibitem

\bibitem{molsch11}
\begin{barticle}[mr]
\bauthor{\bsnm{Molchanov},~\bfnm{Ilya}\binits{I.}} \AND
  \bauthor{\bsnm{Schmutz},~\bfnm{Michael}\binits{M.}}
(\byear{2011}).
\btitle{Exchangeability-type properties of asset prices}.
\bjournal{Adv. in Appl. Probab.}
\bvolume{43}
\bpages{666--687}.
\bid{doi={10.1239/aap/1316792665}, issn={0001-8678}, mr={2858216}}
\bptok{imsref}%
\end{barticle}
\endbibitem

\bibitem{molsch11m}
\begin{bincollection}[author]
\bauthor{\bsnm{Molchanov},~\bfnm{I.}\binits{I.}} \AND
  \bauthor{\bsnm{Schmutz},~\bfnm{M.}\binits{M.}}
(\byear{2013}).
\btitle{Multiasset derivatives and joint distributions of asset prices}.
In \bbooktitle{Musiela Festschrift}
(\beditor{\bfnm{Yu.}\binits{Y.}~\bsnm{Kabanov}}, ed.).
\blocation{Berlin}: \bpublisher{Springer}.
\bnote{To appear}.
\bptok{imsref}%
\end{bincollection}
\endbibitem

\bibitem{mos02}
\begin{bbook}[mr]
\bauthor{\bsnm{Mosler},~\bfnm{Karl}\binits{K.}}
(\byear{2002}).
\btitle{Multivariate Dispersion, Central Regions and Depth:
The Lift Zonoid Approach}.
\bseries{Lecture Notes in Statistics}
\bvolume{165}.
\blocation{Berlin}: \bpublisher{Springer}.
\bid{doi={10.1007/978-1-4613-0045-8}, mr={1913862}}
\bptok{imsref}%
\end{bbook}
\endbibitem

\bibitem{ros95}
\begin{barticle}[mr]
\bauthor{\bsnm{Rosi{\'n}ski},~\bfnm{Jan}\binits{J.}}
(\byear{1995}).
\btitle{On the structure of stationary stable processes}.
\bjournal{Ann. Probab.}
\bvolume{23}
\bpages{1163--1187}.
\bid{issn={0091-1798}, mr={1349166}}
\bptok{imsref}%
\end{barticle}
\endbibitem

\bibitem{ros00}
\begin{barticle}[mr]
\bauthor{\bsnm{Rosi{\'n}ski},~\bfnm{Jan}\binits{J.}}
(\byear{2000}).
\btitle{Decomposition of stationary {$\alpha $}-stable random fields}.
\bjournal{Ann. Probab.}
\bvolume{28}
\bpages{1797--1813}.
\bid{doi={10.1214/aop/1019160508}, issn={0091-1798}, mr={1813849}}
\bptok{imsref}%
\end{barticle}
\endbibitem

\bibitem{ros06}
\begin{barticle}[mr]
\bauthor{\bsnm{Rosi{\'n}ski},~\bfnm{Jan}\binits{J.}}
(\byear{2006}).
\btitle{Minimal integral representations of stable processes}.
\bjournal{Probab. Math. Statist.}
\bvolume{26}
\bpages{121--142}.
\bid{issn={0208-4147}, mr={2301892}}
\bptok{imsref}%
\end{barticle}
\endbibitem

\bibitem{samtaq94}
\begin{bbook}[mr]
\bauthor{\bsnm{Samorodnitsky},~\bfnm{Gennady}\binits{G.}} \AND
  \bauthor{\bsnm{Taqqu},~\bfnm{Murad~S.}\binits{M.S.}}
(\byear{1994}).
\btitle{Stable Non-{G}aussian Random Processes:
Stochastic Models with Infinite Variance}.
\bseries{Stochastic Modeling}.
\blocation{New York}: \bpublisher{Chapman \& Hall}.
\bid{mr={1280932}}
\bptok{imsref}%
\end{bbook}
\endbibitem



\bibitem{sat99}
\begin{bbook}[mr]
\bauthor{\bsnm{Sato},~\bfnm{Ken-iti}\binits{K.i.}}
(\byear{1999}).
\btitle{L\'evy Processes and Infinitely Divisible Distributions}.
\bseries{Cambridge Studies in Advanced Mathematics}
\bvolume{68}.
\blocation{Cambridge}: \bpublisher{Cambridge Univ. Press}.
\bid{mr={1739520}}
\bptok{imsref}%
\end{bbook}
\endbibitem

\bibitem{sat00}
\begin{bmisc}[author]
\bauthor{\bsnm{Sato},~\bfnm{Ken-iti}\binits{K.i.}}
(\byear{2000}).
\bhowpublished{Density transformation in {L\'evy} processes. MaPhySto, Lecture
  Notes, 7}.
\bptok{imsref}%
\end{bmisc}
\endbibitem

\bibitem{schm10qf}
\begin{barticle}[mr]
\bauthor{\bsnm{Schmutz},~\bfnm{Michael}\binits{M.}}
(\byear{2011}).
\btitle{Semi-static hedging for certain {M}argrabe-type options with barriers}.
\bjournal{Quant. Finance}
\bvolume{11}
\bpages{979--986}.
\bid{doi={10.1080/14697688.2010.497494}, issn={1469-7688}, mr={2813831}}
\bptok{imsref}%
\end{barticle}
\endbibitem

\bibitem{schn}
\begin{bbook}[mr]
\bauthor{\bsnm{Schneider},~\bfnm{Rolf}\binits{R.}}
(\byear{1993}).
\btitle{Convex Bodies: The {B}runn--{M}inkowski Theory}.
\bseries{Encyclopedia of Mathematics and Its Applications}
\bvolume{44}.
\blocation{Cambridge}: \bpublisher{Cambridge Univ. Press}.
\bid{doi={10.1017/CBO9780511526282}, mr={1216521}}
\bptok{imsref}%
\end{bbook}
\endbibitem

\bibitem{teh09}
\begin{barticle}[mr]
\bauthor{\bsnm{Tehranchi},~\bfnm{Michael~R.}\binits{M.R.}}
(\byear{2009}).
\btitle{Symmetric martingales and symmetric smiles}.
\bjournal{Stochastic Process. Appl.}
\bvolume{119}
\bpages{3785--3797}.
\bid{doi={10.1016/j.spa.2009.07.007}, issn={0304-4149}, mr={2568296}}
\bptok{imsref}%
\end{barticle}
\endbibitem

\bibitem{wanstoev10}
\begin{barticle}[mr]
\bauthor{\bsnm{Wang},~\bfnm{Yizao}\binits{Y.}} \AND
  \bauthor{\bsnm{Stoev},~\bfnm{Stilian~A.}\binits{S.A.}}
(\byear{2010}).
\btitle{On the association of sum- and max-stable processes}.
\bjournal{Statist. Probab. Lett.}
\bvolume{80}
\bpages{480--488}.
\bid{doi={10.1016/j.spl.2009.12.001}, issn={0167-7152}, mr={2593589}}
\bptok{imsref}%
\end{barticle}
\endbibitem

\bibitem{wanstoev10s}
\begin{barticle}[mr]
\bauthor{\bsnm{Wang},~\bfnm{Yizao}\binits{Y.}} \AND
  \bauthor{\bsnm{Stoev},~\bfnm{Stilian~A.}\binits{S.A.}}
(\byear{2010}).
\btitle{On the structure and representations of max-stable processes}.
\bjournal{Adv. in Appl. Probab.}
\bvolume{42}
\bpages{855--877}.
\bid{doi={10.1239/aap/1282924066}, issn={0001-8678}, mr={2779562}}
\bptok{imsref}%
\end{barticle}
\endbibitem

\end{thebibliography}
\end{document}